\documentclass[12pt]{article}%
\usepackage{multicol,multirow}
\usepackage{amsmath,amssymb,amsfonts}
\usepackage{mathrsfs}
\usepackage{rotating}
\usepackage[numbers,sort&compress]{natbib}
\usepackage{amsthm}
\usepackage{lipsum}
\newtheorem{theorem}{Theorem}[section]
\newtheorem{lemma}[theorem]{Lemma}
\newtheorem{corollary}[theorem]{Corollary}
\newtheorem{proposition}[theorem]{Proposition}

\newtheorem{definition}[theorem]{Definition}
\newtheorem{remark}[theorem]{Remark}
\newtheorem{example}[theorem]{Example}

\numberwithin{equation}{section}
\usepackage{amstext}
\usepackage{latexsym}
\usepackage{amssymb}
\usepackage{amsmath,amsthm,centernot}
\usepackage[hidelinks]{hyperref}
\usepackage{ifpdf}
\usepackage[pagewise]{lineno}
\allowdisplaybreaks

\newcommand{\R}{\mathbb{R}}
\newcommand{\C}{\mathbb{C}}
\newcommand{\N}{\mathbb{N}}
\newcommand{\Z}{\mathbb{Z}}
\newcommand{\LL}{\mathcal{L}}
\newcommand{\eps}{\varepsilon}
\newcommand{\re}{\operatorname{Re}}
\newcommand{\im}{\operatorname{Im}}
\newcommand{\ran}{\operatorname{Ran}}
\newcommand{\ap}{\mathrm{ap}}

\newcommand{\tn}[1]{\left|\!\left|\!\left|#1\right|\!\right|\!\right|}

\begin{document}
	
		\title{Positive expansiveness of operator semigroups: orbit growth and spectral criteria}
	\author{Huynh Minh Hien\\		
		{E-mail: huynhminhhien@qnu.edu.vn},\\
		Department of Mathematics and Statistics, Quy Nhon University,\\ {170 An Duong Vuong}, {Quy Nhon Nam}, Gia Lai,  {Vietnam}\\
	} \date{}
	
		\maketitle



\begin{abstract}
	This paper studies positive expansiveness and positive kinematic expansiveness for strongly continuous semigroups of bounded linear operators on complex Banach spaces. We prove that positive kinematic expansiveness is equivalent to the unboundedness of every nonzero orbit, while positive expansiveness additionally requires the generator to be bounded from below. We also show that continuous changes of time can be replaced by increasing homeomorphisms without changing positive expansiveness. Uniform exponential expansion implies positive expansiveness without any surjectivity assumption; we characterize it through the minimum modulus, the lower growth bound, and the approximate point spectrum of a single semigroup operator. Under a uniform positive lower bound on the semigroup, positive expansiveness is equivalent to positive kinematic expansiveness together with boundedness below of the generator. Diagonal semigroups on $\ell_p$ then show that uniform exponential expansion, positive expansiveness and its kinematic version form a strictly decreasing hierarchy. Finally, we derive necessary spectral conditions for bounded generators and show that all three properties coincide in finite dimensions.
\end{abstract}
	
	\maketitle
	
	\tableofcontents
 
\section{Introduction}

The main goal of the study of dynamical systems is to understand the long time behavior of orbits, and expansiveness has proved to be a powerful tool for studying the chaotic behavior of dynamical systems. 
The study of expansive homeomorphisms started in 1950 with the work of Utz \cite{Utz}. 
In 1972, Bowen and Walters \cite{BW} extended the definition to flows on compact metric spaces. Since two points on the same orbit trace out the same curve at slightly different times, they compared trajectories up to a continuous function on $\R$ vanishing at the origin, the so-called reparametrization. They also showed that, for flows without fixed points, these continuous functions can be replaced by increasing homeomorphisms.
Since then, there have been various types of expansiveness for flows, depending on the reparametrization, including Artigue's kinematic expansiveness \cite{Art}, which compares trajectories at equal times, and the rescaled  expansiveness of Wen and Wen \cite{WW}, which adapts the threshold to the speed of the flow; see \cite{JNY} for the spectral decomposition of rescaling expansive flows, \cite{LMV} for asymptotic expansivity, and the monograph \cite{LMN2} for measure-theoretic versions. For semiflows the situation is different. Orbits can merge in forward time, and the Bowen--Walters definition becomes extremely restrictive: Herrero and Jaque \cite{HJ} showed that a semiflow on a metric space that is expansive in the sense of Alves, Carvalho and Siqueira is trivial, and that a positively expansive semiflow on a compact metric space consists of finitely many closed orbits, unbranched tails and isolated singularities. Lee and Rojas \cite{LR} introduced eventual expansiveness, in which nearby orbits are only required to coincide eventually, and proved a similar finiteness result for injective eventually expansive semiflows on compact spaces. Thus, on compact spaces, expansive semiflows carry little dynamics. 

	Linear operators offer a very different picture. The phase space is no longer compact, separation can be measured relative to the size of the initial point, and expansiveness becomes a property of individual orbits governed by spectral theory. An automorphism of a finite-dimensional space is expansive if and only if it has no eigenvalue of modulus one \cite{Eis}. Eisenberg and Hedlund \cite{EH} related expansiveness of automorphism on Banach spaces to the spectrum and the approximate point spectrum, and Hedlund \cite{Hed} proved that an automorphism is uniformly expansive if and only if its approximate point spectrum misses the unit circle. Without invertibility, Bernardes et al. \cite{BCDMP} showed that an operator is positively expansive if and only if every nonzero forward orbit is unbounded, and related uniform positive expansiveness to the position of the approximate point spectrum relative to the closed unit disc. These ideas have since been connected with shadowing, hyperbolicity and structural stability \cite{BM,CGP,Mai,LMN}, and with Cauchy problems on Banach spaces \cite{YSMN}. In discrete time, expansiveness of linear maps is thus largely understood through two objects: the growth of orbits and the approximate point spectrum.

A semigroup $T=(T(t))_{t\ge0}$ on a Banach space $X$, with generator $A$, is simultaneously a semiflow on a noncompact space and a continuous family of linear operators, and it is natural to ask how much of the discrete theory survives when time becomes continuous and reparametrizations may be changed. Lee and Morales \cite{LM} took the first step by introducing positive expansiveness in the sense of Bowen--Walters for $C_0$-semigroups, with continuous reparametrization and a separation threshold proportional to $\|x\|$. They proved that the point spectrum of the generator of a positively expansive semigroup lies in the open right half-plane, and that such a semigroup is neither uniformly bounded, nor equicontinuous, nor chaotic. As a sufficient condition, they showed that the resolvent estimate $\|(\gamma I-A)^{-1}\|\le(\alpha-\gamma)^{-1}$, $\gamma<\alpha$, with $\alpha>0$, yields invertible operators with $\|T(t)^{-1}\|\le e^{-\alpha t}$, and hence positive expansiveness, and they applied this to the Lasota equation. Several questions arise naturally. Firstly, their sufficient condition requires every $T(t)$ to be invertible, which excludes natural expanding semigroups whose operators are not surjective, such as weighted right translations. Secondly, their necessary spectral condition involves only the point spectrum, which may be empty, and no characterization in terms of the approximate point spectrum, comparable to \cite{Hed,BCDMP}, is available. Thirdly, the role of reparametrization, and the gap between positive expansiveness and exponential growth, are unclear. Finally, the Lasota semigroup is not strongly continuous on the space used in \cite{LM}. 
	
To address these questions, we study \textit{positive expansiveness} together with its kinematic version, \textit{positive kinematic expansiveness}, a one-side linear analogue of Artigue's notion, in which trajectories are compared at equal times. As a sufficient condition we use \textit{uniform exponential expansion}: there are constants $M,a>0$ such that
\begin{equation}\label{intro}
	\|T(t)x\|\ge M e^{at}\|x\|\quad\mbox{for}\quad x\in X,\ t\ge0.
\end{equation}
The corresponding estimate with $a=0$ is a \textit{uniform lower bound}
\begin{equation}\label{ulb} \|T(t)x\|\ge M\|x\| \quad\mbox{for}\quad x\in X,\ t\ge0.
\end{equation}
These three properties satisfy
\begin{equation}\label{3notion}
	\begin{aligned}
		&\mbox{uniform exponential expansion}\ \Longrightarrow\ \mbox{positive expansiveness}\\
		&\qquad\quad\Longrightarrow\ \mbox{positive kinematic expansiveness},
	\end{aligned}
\end{equation}
and neither implication can be reversed. Our results can be summarized in three points.

\textit{Orbits and the generator.} On complex spaces, positive kinematic expansiveness holds if and only if every nonzero orbit is unbounded, and the necessary conditions of \cite{LM} already hold for it. Allowing changes of time adds a further necessary condition, namely that the generator be bounded from below, and under the uniform lower bound \eqref{ulb} these two conditions characterize positive expansiveness.

\textit{Changes of time.} Positive expansiveness is unchanged if the continuous reparametrizations in its definition are restricted to increasing homeomorphisms. when the flow has no singularities, and more generally when its set of singularities is open \cite{Hien}. For $C_0$-semigroups, by contrast, the origin is always a fixed point that is not isolated, and its role is taken over by the condition that the generator be bounded from below. In addition, the space is no longer compact.

\textit{Uniform exponential expansion.} Condition \eqref{intro} implies positive expansiveness even when no operator $T(t)$ is surjective, and it can be checked on a single operator: it holds if and only if the approximate point spectrum of some $T(\tau)$ lies outside the closed unit disc. In the invertible case the single-time bound $\|T(\tau)^{-1}\|<1$ is sufficient for positive expansiveness, so the hypothesis of Lee-Morales \cite{LM} is a special case. Moreover, we introduce the continuous time counterpart of the lower spectral radius of an operator, which is called the \textit{lower growth bound} of a semigroup. This yields criteria for uniform exponential expansion, and hence for positive expansiveness, in terms of a single operator of the semigroup or of the generator. 
 Unlike positive expansiveness, \eqref{intro} persists under small bounded perturbations of the generator.

The paper is organized as follows. In Section~\ref{sec2} we introduce the two notions of expansiveness and the quantities that measure expansion from below, namely the minimum modulus, the lower spectral radius and the lower growth bound, and we express the condition $0\notin\sigma_\ap(A)$ as a uniform lower bound on short time displacements. Section~\ref{sec3} is concerned with necessary conditions. We first show that positive kinematic expansiveness depends only on the growth of individual orbits; we then prove that positive expansiveness forces the generator to be bounded from below, and that continuous reparametrizations may be replaced by increasing homeomorphisms without changing the notion. In Section~\ref{sec4} we turn to uniform exponential expansion. We prove that it implies positive expansiveness without any surjectivity assumption and that it can be read off from a single operator of the semigroup; we also treat the invertible case, which contains the sufficient condition of Lee and Morales, give a Lyapunov-type characterization on Hilbert spaces, and show that uniform exponential expansion survives small bounded perturbations. Section~\ref{sec5} deals with semigroups that only satisfy a uniform lower bound. Here the necessary conditions of Section~\ref{sec3} turn out to be sufficient, and we apply this characterization to diagonal semigroups, which show that both implications in \eqref{3notion} are strict. Finally, in Section~\ref{sec6} we give criteria in terms of the spectrum of the generator; for bounded generators we compare uniform exponential expansion with positive expansiveness and show that the gap between them closes in finite dimensions but not in general.

\section{Preliminaries}\label{sec2}
In this section we fix notation, recall the definition of positive expansiveness and its kinematic version, and introduce the quantities used to measure expansion from below: the minimum modulus, the lower spectral radius and the lower growth bound. Their relation to the approximate point spectrum provides the spectral tools for the subsequent sections.

\subsection{Semigroups and generators}\label{subsec21}
Throughout this paper, $X$ denotes a nonzero complex Banach space and $\LL(X)$ denotes the space of bounded linear operators on $X$. A \textit{semigroup} is a family 
$T=(T(t))_{t\ge0}\subset\LL(X)$ with
$$T(0)=I\quad\mbox{and}\quad T(t+s)=T(t)T(s)\quad\mbox{for}\quad t,s\ge0;$$ it is a \textit{$C_0$-semigroup}, or \textit{strongly continuous},  if $T(t)x\to x$ as $t\to 0^+$ for every $x\in X$, and it is called \textit{uniformly continuous} if $\|T(t)-I\|\to0$ as $t\to 0^+$.  Its generator $A$ is given by
\[
D(A)=\Big\{x\in X:\lim_{t\to 0^+}\frac{T(t)x-x}{t}\ \mbox{exists}\Big\},\quad Ax=\lim_{t\to 0^+}\frac{T(t)x-x}{t}.
\]
Then $A$ is linear, closed on $D(A)$, which is dense in $X$. We write $\sigma(A)$, $\sigma_p(A)$, $\sigma_\ap(A)$ and $\rho(A)$ for the spectrum, point spectrum, approximate point spectrum and resolvent set of $A$, respectively. A $C_0$-semigroup is bounded on compact time intervals and has continuous orbits; it is uniformly continuous if and only if its generator is bounded, in which case $T(t)=e^{tA}$. We refer to \cite{EN,Pazy} for these and other standard facts.

The next lemma expresses the condition $0\notin\sigma_\ap(A)$, i.e. boundedness from below of the generator, as a uniform lower bound on the displacement of the semigroup over short times. It is in this form that the condition enters the proofs of Theorems~\ref{mthm0} and~\ref{mainthm2}.
\begin{lemma}\label{eqlm}
	Let $T$ be a $C_0$-semigroup with generator $A$. The following assertions are equivalent:
	\begin{enumerate}
		\item[(i)] $0\notin\sigma_\ap(A)$;
		\item[(ii)] there is $\kappa>0$ such that $\|Ax\|\ge\kappa\|x\|$ for $x\in D(A)$;
		\item[(iii)] for every $\eps>0$ there is $\kappa_\eps>0$ such that
		\begin{equation}\label{cleqsup}
			\sup_{0\le t\le\eps}\|T(t)x-x\|\ \ge\ \kappa_\eps\|x\|\quad\mbox{for}\quad x\in X.
		\end{equation}
	\end{enumerate}
\end{lemma}

\begin{proof}
	(i)~$\Rightarrow$~(ii). Recall from \cite[Lemma VI.1.9]{EN} that \begin{equation*}
		\sigma_\ap(A)=\{\lambda\in\C:\ \exists\, x_n\in D(A),\ \|x_n\|=1,\ \|(A-\lambda I)x_n\|\to 0\}.
	\end{equation*}
	Therefore, $0\in\sigma_{\mathrm{ap}}(A)$ if and only if 
	there exist $x_n\in D(A)$ with $\|x_n\|=1$ and $\|Ax_n\|\to 0$.
	Suppose on the contrary that there were no $\kappa>0$ satisfying
	\[
	\|Ax\|\ge\kappa\|x\|\quad\mbox{for}\quad x\in D(A).
	\]
	Then, for each $n\ge1$, we could choose $y_n\in D(A)\setminus\{0\}$ such that $\|Ay_n\|<\frac1n\|y_n\|$.
	Setting $x_n=\frac{y_n}{\|y_n\|}$, we obtain
	\[
	\|x_n\|=1\quad\mbox{and}\quad 
	\|Ax_n\|<\frac1n\to0\quad\mbox{as}\quad n\to\infty.
	\]
	This would imply $0\in\sigma_{\mathrm{ap}}(A)$, a contradiction.
	
	(ii)~$\Rightarrow$~(i). Suppose  that there is $\kappa>0$ such that $\|Ax\|\ge \kappa\|x\|$ for $x\in D(A)$. Then, every unit vector $x\in D(A)$ satisfies $\|Ax\|\ge\kappa$, so no sequence of unit vectors can satisfy $Ax_n\to 0$. Hence $0\notin\sigma_{\mathrm{ap}}(A)$.
	
	(ii)~$\Rightarrow$~(iii). Suppose \eqref{cleqsup} fails for some $\eps>0$ and every $\kappa_\eps>0$. Then there are unit vectors $z_n$ with $\eta_n:=\displaystyle\sup_{0\le t\le\eps}\|T(t)z_n-z_n\|\to 0$. Put $w_n=\frac1\eps\int_0^\eps T(t)z_n\,dt$. By \cite[Lemma II.1.3]{EN}, $w_n\in D(A)$ and $Aw_n=\frac1\eps(T(\eps)z_n-z_n)$, so \begin{equation}\label{eq71} 
		\|Aw_n\|\le\frac{\eta_n}\eps\to 0\quad\mbox{as}\quad n\to\infty.
	\end{equation} Moreover, \[\|w_n-z_n\|\le\frac1\eps\int_0^\eps\|T(t)z_n-z_n\|dt\le\eta_n,\] so 
	\[\|w_n\|\ge1-\eta_n\to1\quad\mbox{as}\quad n\to\infty.\]
	Together with \eqref{eq71}, it follows that
	\begin{equation}
		\left\|A\left(\frac{w_n}{\|w_n\|}\right)\right\|=\frac{\|A w_n\|}{\|w_n\|}\to 0\quad\mbox{as}\quad n\to\infty,
	\end{equation} i.e.\ $0\in\sigma_\ap(A)$, contradicting (ii).
	
	(iii)~$\Rightarrow$~(i). Fix $\eps>0$ and $\kappa_\eps$ satisfying \eqref{cleqsup}. Suppose $0\in\sigma_\ap(A)$, and let $x_n\in D(A)$ be unit vectors with $Ax_n\to 0$ as $n\to\infty$. With $K_\eps=\displaystyle\sup_{0\le t\le\eps}\|T(t)\|$, using $\frac{d}{dt}T(t)x_n=T(t)Ax_n$, we have, for $0\le t\le\eps$, \[\|T(t)x_n-x_n\|=\left\|\int_0^tT(u)Ax_n\,du\right\|\le |t|K_\eps\|Ax_n\|\le \eps K_\eps\|Ax_n\|\to 0,\]
	as $n\to\infty$, which contradicts \eqref{cleqsup}.		
\end{proof}

\subsection{Positive expansiveness and its kinematic version}\label{subsec22}
Denote by ${\cal C}$ the set of continuous functions $s:[0,\infty)\to [0,\infty)$ with $s(0)=0$.

\begin{definition}[\cite{LM}]\label{LMdef} \rm 
	A semigroup $T$ is called \textit{positively expansive} if for every $\eps>0$ there exists $\delta>0$ such that, whenever $x,y\in X$ and $s\in {\cal C}$ satisfy
	\begin{equation}\label{LMdefeq}
		\|T(t)x-T(s(t))y\|<\delta\|x\|\quad\mbox{for}\quad t\ge 0,
	\end{equation}
	there is $t_*\in[0,\eps)$ for which $y=T(t_*)x$ or $x=T(t_*)y$.
\end{definition}
In fact, in \cite{LM}, the condition reads $t_*\in[0,\varepsilon]$, but the two definitions are equivalent: the version with $[0,\varepsilon)$ trivially implies the one with $[0,\varepsilon]$, and conversely one applies the closed version with $\varepsilon/2$ in place of $\varepsilon$. We can assume that $\|x\|=1$ and \eqref{LMdefeq} becomes
\[\|T(t)x-T(s(t))y\|<\delta\quad\mbox{for}\quad t\ge 0. \] 
Restricting the reparametrization in Definition~\ref{LMdef} to the identity gives the following kinematic version.
\begin{definition}\label{pkedef}\rm
	A semigroup $T$ is \textit{positively kinematic expansive} if for every $\eps>0$ there exists $\delta>0$ such that, whenever $x,y\in X$ satisfy
	\begin{equation}\label{closeeq}
		\|T(t)x-T(t)y\|<\delta\|x\|\quad\mbox{for}\quad t\ge 0,
	\end{equation}
	there is $t_*\in[0,\eps)$ for which $y=T(t_*)x$ or $x=T(t_*)y$.
\end{definition}
This is  the linear, one-sided counterpart of Artigue's kinematic expansiveness for flows \cite{Art}. It is obvious that 
\[\mbox{positive expansiveness}\Longrightarrow\mbox{positive kinematic expansiveness}.\]

\subsection{Minimum modulus and lower spectral radius}\label{subsec23}
The lower estimates used throughout the paper are expressed through the \textit{minimum modulus}. For $L\in\LL(X)$, this quantity is defined by
\begin{equation}\label{Ldn}
	\tn{L}:=\inf_{\|x\|=1}\|Lx\|.
\end{equation}
It satisfies, for $L,L_1,L_2\in\LL(X)$,
\begin{equation}\label{L1L2}
	\tn{L_1L_2} \ge \tn{L_1}\,\tn{L_2},\quad \tn{L_1L_2} \le \|L_1\|\,\tn{L_2}.
\end{equation}
The first inequality in \eqref{L1L2} follows from $\|L_1L_2x\|\ge\tn{L_1}\|L_2x\|$, the second from $\|L_1L_2x\|\le\|L_1\|\|L_2x\|$.

In the case that $L\in\LL(X)$ is invertible, since $y=Lx$ ranges over $X\setminus\{0\}$ as $x$ does, we obtain
\begin{equation*}
	\|L^{-1}\|
	=\displaystyle\sup_{y\ne0}\frac{\|L^{-1}y\|}{\|y\|}
	=\displaystyle\sup_{x\ne0}\frac{\|x\|}{\|Lx\|}
	=\left(\inf_{\|x\|=1}\|Lx\|\right)^{-1}
	=\tn{L}^{-1}.
\end{equation*}
Consequently 
\begin{equation}\label{L-1}
	\tn{L}=\|L^{-1}\|^{-1}.
\end{equation}

The minimum modulus is closely related to the approximate point spectrum. For $L\in\LL(X)$, the set $\sigma_\ap(L)$ is a nonempty compact set, and 
\begin{equation}\label{sigmaL} \partial \sigma(L)\subset \sigma_\ap(L);
\end{equation}
see \cite[Lemma IV.1.10]{EN} for a proof.
Moreover, 
\begin{equation}\label{Lmeq}
	\tn{L}\le|\mu|\quad\mbox{for}\quad \mu\in\sigma_\ap(L).
\end{equation}
Indeed, if $\|z_n\|=1$ and $(L-\mu I)z_n\to 0$, then $\tn{L}\le\|Lz_n\|\le\|(L-\mu I)z_n\|+|\mu|$, and the assertion follows by letting $n\to\infty$.

To describe the asymptotic behavior of the minimum modulus under iteration, we use the \textit{lower spectral radius} of $L\in\LL(X)$, defined by
\begin{equation}\label{r-}
	r_-(L):=\lim_{n\to\infty}\tn{L^n}^{1/n}.
\end{equation}
To verify that this limit exists, first suppose that $\tn{L}=0$, then $\tn{L^n}\le\|L^{n-1}\|\tn{L}=0$ for all $n\ge1$, so $r_-(L)=0$. If $\tn{L}>0$, then $\tn{L^n}\ge\tn{L}^n>0$, the sequence $\{\log\tn{L^n}\}$ is superadditive by \eqref{L1L2}, and Fekete's lemma gives the existence of the limit, which is finite because $\tn{L^n}\le\|L\|^n$. Since $L-\lambda I$ is bounded below precisely when $\lambda\notin\sigma_\ap(L)$, and $\sigma_\ap(L)$ is closed and nonempty, the radius of the largest open disc centred at $0$ on which $L-\lambda I$ is bounded below equals $\min\{|\lambda|:\lambda\in\sigma_\ap(L)\}$. By \cite[Theorem~3]{MZ}, this radius equals $r_-(L)$, i.e.
\begin{equation}\label{MZeq}
	r_-(L)=\min\{|\lambda|\,:\, \lambda\in\sigma_\ap(L)\}.
\end{equation}
Formula \eqref{MZeq} will connect uniform lower estimates with spectral conditions in Theorems \ref{uepeqthm}~(v) and \ref{geAthm}\,(ii).

\subsection{Lower growth bound}\label{subsec24}
We define the \textit{lower growth bound} of a $C_0$-semigroup $T$ by
\begin{equation}\label{omega-dn}
	\omega_-(T):=\lim_{t\to\infty}\frac1t\log\tn{T(t)}\in[-\infty,\infty),
\end{equation}
with convention $\log 0:=-\infty$. This is the continuous time counterpart of the lower spectral radius. 	
The following lemma is the counterpart, for the minimum modulus, of the familiar properties of the growth bound (see \cite[Proposition~IV.2.2]{EN}); we include the proof for completeness.

\begin{lemma}[Lower growth bound]\label{lgblm}
	Let $T$ be a $C_0$-semigroup. The following hold:
	\begin{enumerate}
		\item[(i)] Either $\tn{T(t)}=0$ for all $t>0$, or $\tn{T(t)}>0$ for all $t\ge 0$.
		\item[(ii)] The limit in \eqref{omega-dn} exists and
		\begin{equation}\label{omega-T}
			\omega_-(T)=\displaystyle\sup_{t>0}\frac1t\log\tn{T(t)}.
		\end{equation}
		\item[(iii)] For every $\tau>0$, $\;r_-(T(\tau))=e^{\tau\omega_-(T)}$.
	\end{enumerate}
\end{lemma}
\begin{proof}
	(i) Suppose $\tn{T(\tau)}>0$ for some $\tau>0$. By the first inequality in \eqref{L1L2}, $\tn{T(n\tau)}\ge\tn{T(\tau)}^n>0$ for every integer $n\ge1$. Given $t>0$, choose $n$ with $n\tau\ge t$. The second inequality in \eqref{L1L2} gives
	\begin{equation*}
		0<\tn{T(n\tau)}=\tn{T(n\tau-t)T(t)}\le\|T(n\tau-t)\|\,\tn{T(t)},
	\end{equation*}
	hence $\tn{T(t)}>0$. For $t=0$ we have $\tn{T(0)}=\tn{I}=1$.
	
	(ii) If $\tn{T(t)}=0$ for all $t>0$, then both sides of \eqref{omega-T} equal $-\infty$. By (i), it remains to consider the case $\tn{T(t)}>0$ for all $t\ge0$. Then
	\[
	f(t):=\log\tn{T(t)},\quad t\ge0,
	\]
	is a real valued function with $f(0)=0$. 
	First, $f$ is superadditive since by the first inequality in \eqref{L1L2},
	\[
	f(t+s)=\log\tn{T(t)T(s)}\ge \log\tn{T(t)}+\log\tn{T(s)} = f(t)+f(s)\quad\mbox{for}\quad t,s\ge0.\] Next, fix $h>0$ and put $K_h:=\sup_{0\le u\le h}\|T(u)\|<\infty$. For $0\le r\le h$, the second inequality in \eqref{L1L2} gives
	\[
	\tn{T(h)}=\tn{T(h-r)T(r)}\le K_h\tn{T(r)},
	\]
	so
	\[
	c_h:=\inf_{0\le r\le h}f(r)\ \ge\ f(h)-\log K_h>-\infty.
	\]
	Fix $h>0$ and write $t=nh+r$ with $n\in\N$ and $0\le r<h$. Superadditivity gives $f(t)\ge nf(h)+c_h$, and since $nh/t\to1$ and $c_h/t\to 0$,
	\begin{equation*}
		\liminf_{t\to\infty}\frac{f(t)}t\ \ge\ \frac{f(h)}h .
	\end{equation*}
	Taking the supremum over $h>0$ yields
	\begin{equation*}
		\liminf_{t\to\infty}\frac{f(t)}t
		\ge\displaystyle\sup_{h>0}\frac{f(h)}h
		\ge\limsup_{t\to\infty}\frac{f(t)}t.
	\end{equation*} Hence the limit exists and equals the supremum. Finally, choose $C\ge1$ and $\omega\in\R$ with $\|T(t)\|\le Ce^{\omega t}$ for $t\ge0$. Then 
	\[f(t)\le\log\|T(t)\|\le\log C+\omega t,\] so the limit is at most $\omega$ and hence it is finite. 
	
	(iii) Fix $\tau>0$. If $\tn{T(t)}=0$ for all $t>0$, then $\tn{T(\tau)^n}=\tn{T(n\tau)}=0$ for every $n\ge1$, so $r_-(T(\tau))=0=e^{\tau\omega_-(T)}$ with the convention $e^{-\infty}=0$. Otherwise, $T(\tau)^n=T(n\tau)$, and the limit in (ii), taken along the sequence $t=n\tau\to\infty$, gives
	\[
	\log r_-(T(\tau))=\lim_{n\to\infty}\frac1n\log\tn{T(n\tau)}=\tau\lim_{n\to\infty}\frac{f(n\tau)}{n\tau}=\tau\,\omega_-(T).
	\]
	Hence $r_-(T(\tau))=e^{\tau\omega_-(T)}$.
\end{proof}
	\section{Necessary conditions and the role of the reparametrization}  \label{sec3}
	
	In this section we determine what expansiveness imposes on a semigroup and its generator. We first treat the kinematic notion and show that it is governed entirely by the growth of individual orbits. We then show that positive expansiveness adds a condition on the generator, namely that it be bounded from below, and that the class of reparametrizations can be reduced to increasing homeomorphisms without changing the notion. The two necessary conditions obtained here are the ingredients of the characterization in Section~\ref{sec5}.
	
	\subsection{Positive kinematic expansiveness}
			
	\begin{lemma}\label{bolm}
		Let $T$ be a $C_0$-semigroup on a Banach space $X$. Suppose that $T$ is positively kinematic expansive and  $x\ne0$ such that $\displaystyle\sup_{t\ge 0}\|T(t)x\|<\infty$. Then for every $\eps>0$ there exist $r\in(0,\eps)$ and $\alpha\in(0,\infty)\setminus\{1\}$ with $T(r)x=\alpha x$.
	\end{lemma}
	
	\begin{proof} Let $x\ne 0$ and $K=\displaystyle\sup_{t\ge 0}\|T(t)x\|<\infty$, For every $\eps$, let $\delta>0$ correspond to $\eps$ in Definition \ref{pkedef}. Choose $c>0$ with $cK<\delta\|x\|$ and put $y=(1+c)x$. Then \[\|T(t)x-T(t)y\|=c\|T(t)x\|\le cK<\delta\|x\| \quad\mbox{for all}\quad t\ge 0.\] 
		Hence $y=T(r)x$ or $x=T(r)y$ for some $r\in[0,\eps)$. Since $y\ne x$, necessarily $r>0$, and $T(r)x=\alpha x$ with $\alpha=1+c$ or $\alpha=(1+c)^{-1}$. It is obvious that $\alpha>0$ and $\alpha\ne 1$.
	\end{proof}

	\begin{theorem}\label{nethm}
		Let $T$ be a $C_0$-semigroup with generator $A$. If $T$ is positively kinematic expansive, then:
		\begin{enumerate}
			\item[(i)] $T(t)$ is injective for every $t>0$;
			\item[(ii)] $T$ has no nonzero periodic points; in particular $\ker A=\{0\}$;
			\item[(iii)] $\sigma_p(A)\subset\{\lambda\in\C:\re\lambda>0\}$.
		\end{enumerate}
	\end{theorem}

Part (iii) was proved by Lee and Morales \cite[Theorem 1.3(1)]{LM} for positively expansive semigroups; the same argument works for the kinematic notion. We include a proof because the argument is reused in Section \ref{sec5}, and because it differs slightly from theirs.
	
	\begin{proof}
		(i) Suppose $T(t_0)x=0$ for some $t_0>0$ and $x\neq 0$. Then $T(t)x=0$ for $t\ge t_0$, so the orbit is bounded. Lemma \ref{bolm} gives $r>0$ and $\alpha>0$ with $T(r)x=\alpha x$, whence $T(nr)x=\alpha^nx\ne0$ for all $n$. Taking $nr\ge t_0$ gives a contradiction.
		
		(ii) Suppose $T(p)x=x$ with $p>0$ and $x\ne0$. The orbit is $p$-periodic, hence bounded. If $T(u)x=0$ for some $u>0$, choose $n$ with $np\ge u$; then 
		\[x=T(np)x=T(np-u)T(u)x=0.\] It follows that
		$T(t)x\ne0$ for all $t\ge 0$. The orbit is continuous and nonvanishing on $[0,p]$ and $p$-periodic, so it has a positive minimum:
		\[0<m\le\|T(t)x\|\le K<\infty\quad\mbox{for all}\quad t\ge 0.\]
		with some $m,K>0$.  Lemma \ref{bolm} gives $T(r)x=\alpha x$ with $\alpha\ne1$, so $\|T(nr)x\|=\alpha^n\|x\|$ tends to $\infty$ or to $0$ as $n\to\infty$, contradicting the bounds $m$ and $K$. Therefore, $T(p)x\ne x$ for $p>0$ and $x\ne 0$. 
		
		Finally, let $x\in\ker A$. For $x\in D(A)$ the orbit $t\mapsto T(t)x$ is continuously differentiable with derivative $T(t)Ax$ (see \cite[Lemma II.1.3]{EN}). Hence $\frac{d}{dt}T(t)x=T(t)Ax=0$, so $T(t)x=x$ for all $t\ge 0$; in particular $T(1)x=x$, and $x=0$ by what was just shown. So, ${\rm ker} A=\{0\}$. 
		
		(iii) Let $\lambda\in\sigma_p(A)$ and $Ax=\lambda x$ with $x\ne0$. Suppose on the contrary that $\alpha:=\re\lambda\le0$; write $\beta=\im\lambda$. Then $T(t)x=e^{\lambda t}x$. Fix $\eps=1$ and let $\delta=\delta(\eps)\in(0,1)$ as in Definition \ref{pkedef}. Choose $c\in(0,\delta)$ so small that $c<\pi$ and, if $\alpha<0$ and $\beta\neq0$, $c<|\alpha|/|\beta|$. For $y=(1+ic)x$
		\begin{equation*}
		\|T(t)x-T(t)y\|=c\,e^{\alpha t}\|x\|\le c\|x\|<\delta\|x\|\quad \mbox{for all}\quad t\ge 0,
		\end{equation*}
		so there is $r\in[0,1)$ with $y=T(r)x$ or $x=T(r)y$. 
		
		If $y=T(r)x$, then $e^{\lambda r}=1+ic$ implies $e^{\alpha r}=\sqrt{1+c^2}>1$, which is impossible for $\alpha\le0$ and $r\ge 0$. 
		
		If $x=T(r)y$, then $e^{\lambda r}(1+ic)=1$ yields $e^{\alpha r}=(1+c^2)^{-1/2}<1$, so $\alpha<0$ and $0<r=\frac{\log(1+c^2)}{2|\alpha|}\le\frac{c^2}{2|\alpha|}$. Taking arguments, $r\beta+\arctan c\in2\pi\Z$. But $0<\arctan c\le c<\pi$ and $|r\beta|\le|\beta|c^2/(2|\alpha|)<c/2\le\arctan c$, using  $\arctan c\ge c/2$ for $0<c\le1$, so $r\beta+\arctan c\in(0,2\pi)$, a contradiction. Hence, $\re\lambda>0$, which was to be shown.
	\end{proof}
	The next result characterizes positive kinematic expansiveness.
	\begin{theorem}\label{ubothm}
		Let $T$ be a $C_0$-semigroup in a complex Banach space $X$. Then $T$ is positively kinematic expansive if and only if every nonzero orbit is unbounded, i.e. $\displaystyle\sup_{t\ge 0}\|T(t)x\|=\infty$ for every $x\ne0$.
	\end{theorem}
	
	\begin{proof} \textit{Sufficiency.} Suppose that every nonzero orbit of $T$ is unbounded. For $x,y\in X$ and 
		\[\|T(t)x-T(t)y\|<\delta\|x\|\quad\mbox{for}\quad t\ge0.\]
		Since $T(t)$ is linear, 
		\[\|T(t)(x-y)\|<\delta \|x\|\quad\mbox{for}\quad t\ge0. \]
		This means that the orbit of $x-y$ is bounded and hence $x=y$, i.e. $y=T(t_*)x$ with $t_*=0$. 
		
		\textit{Necessity}. Suppose $x\ne 0$ and $K:=\displaystyle\sup_{t\ge 0}\|T(t)x\|<\infty$. We claim that   \begin{equation}\label{Ttx}
			T(t)x\in \C x:=\{\lambda x\,:\, \lambda\in\C\}\quad\mbox{for all}\quad t\ge 0.
		\end{equation}
		Define \[S=\{t\ge 0:\ T(t)x\in\C x\}\] and see that \eqref{Ttx} holds if $S=[0,\infty)$.  We first verify that $S$ contains arbitrarily small positive numbers. Let $\eps>0$, let $\delta$ correspond to $\eps$ in Definition \ref{pkedef}, and choose $0<c<1$ satisfying $cK<\delta\|x\|$. For $y=(1+c)x$
		\[ \|T(t)x-T(t)y\|=c\|T(t)x\|<\delta\|x\|\quad\mbox{for all}\quad  t\ge 0.\] 
		Since $T$ is positively kinematic expansive, there is $r\in [0,\eps)$ such that $y=T(r)x$ or $x=T(r)y$. Since $y\ne x$, $r>0$, and in either case $T(r)x\in\{(1+c)x,(1+c)^{-1}x\}\subset\C x$. Thus $S\cap(0,\eps)\ne\varnothing$ for every $\eps>0$.
		
		Next, observe that $S$ is closed, because $t\mapsto T(t)x$ is continuous due to the strong continuity of $T$ and $\C x$ is closed. In addition, $S$ is closed under addition. Indeed,  if $t,t'\in S$, then $T(t)x=\gamma x$ and $T(t')x=\gamma'x$ for some $\gamma,\gamma'\in \C$. Using the semigroup property and linearity,
		\[ 
		T(t+t')x
		=T(t)T(t')x
		=T(t)(\gamma'x)
		=\gamma'T(t)x
		=\gamma'\gamma x
		\] and hence $t+t'\in S$.
		
		Finally, take any $r\in [0,\infty)$, choose $\{r_k\}_{k\ge 1}\subset S$ with $r_k\to 0^+$ as $k\to\infty$. Such a sequence $\{r_k\}$ does exist since $S\cap (0,\eps)\ne\varnothing$ as shown above.  Let $n_k=[ r/r_k]$; then $n_kr_k\in S$ and $n_kr_k\to r$ as $k\to\infty$, so $r\in S$ since $S$ is closed. Hence $S=[0,\infty)$.
		
		We are in a position to show that $K<\infty$ is impossible. Property \eqref{Ttx} shows that for every $t\ge 0$, there is $\gamma(t)\in \C$ such that $T(t)x=\gamma(t)x$. Define  function $\gamma: [0,\infty)\to \C, t\mapsto \gamma(t)$ and choose a bounded functional $\varphi\in X^*$ with $\varphi(x)=1$. Then $\gamma(t)=\varphi(T(t)x)$, so $\gamma$ is continuous due to the strong continuity of $T$. 
		Moreover, \[\gamma(t+s)=\gamma(t)\gamma(s)\quad\mbox{for}\quad  t,s\in [0,\infty)\]
		and  $\gamma(t)\ne0$ for every $t\ge 0$ because $T(t)$ is injective by Theorem \ref{nethm}\,(i). Hence $\gamma(t)=e^{\lambda t}$ for some $\lambda\in\C$. Then 
		\[ \lim_{t\to 0^+}\frac{T(t)x-x}{t} =\lambda x,\] so $x\in D(A)$ and $Ax=\lambda x$, where $A$ is the generator of $T$. It follows that $\lambda\in\sigma_{p}(A)$. By Theorem \ref{nethm}\,(iii), $\re\lambda>0$, and then $\|T(t)x\|=e^{t\re\lambda}\|x\|\to\infty$, contradicting $K<\infty$.

	\end{proof}
	
	The following property is established in \cite{LM} for positively expansive semigroups.
	\begin{corollary}\label{kinco}
		A positively kinematic expansive $C_0$-semigroup $T$ on a complex Banach space $X\ne\{0\}$ is neither uniformly bounded, nor equicontinuous, nor chaotic.
	\end{corollary}
	
	\begin{proof}
		If $\displaystyle\sup_{t\ge 0}\|T(t)\|<\infty$, every orbit is bounded, contradicting Theorem \ref{ubothm}. Equicontinuity is equivalent to uniform boundedness \cite[Lemma 2.3]{LM}. Finally, $T$ is chaotic if it has a dense orbit and a dense set of periodic points \cite{GEM}; since $X\ne\{0\}$, a dense set of periodic points contains a nonzero periodic point, which is excluded by Theorem \ref{nethm}\,(ii).
	\end{proof}

	\begin{remark}\rm
		Theorem~\ref{ubothm} is considerably stronger than the failure of uniform boundedness in Corollary~\ref{kinco}. By the uniform boundedness principle, $\sup_{t\ge0}\|T(t)\|=\infty$ only means that some orbit is unbounded, whereas Theorem~\ref{ubothm} requires every nonzero orbit to be unbounded. For instance, the semigroup $T$ on $\C^2$ given by
		\[
		T(t)(x_1,x_2)=(e^{t}x_1,\;x_2),\quad (x_1,x_2)\in\C^2,\ t\ge0,
		\]
		is not uniformly bounded, since $\|T(t)\|=e^{t}$, but the orbit of $(0,1)$ is constant, so $T$ is not positively kinematic expansive;  this is detected by Theorem~\ref{ubothm} but not by Corollary~\ref{kinco}.  {\hfill$\Diamond$}
	\end{remark}
	
	\subsection{Positive expansiveness and the generator}
	The next necessary condition concerns the generator as a whole. The parametrization $s$ in the proof is not the identity (see Definition \ref{LMdef}).
	
	\begin{theorem}\label{nAxthm} Let $T$ be a $C_0$-semigroup on a Banach space $X$ with generator $A$. If $T$ is positively expansive, then $0\notin\sigma_\ap(A)$.
	\end{theorem}
	
	\begin{proof}
		Suppose $0\in\sigma_\ap(A)$ and choose unit vectors $z_n\in D(A)$ with $Az_n\to 0$. Let $K=\displaystyle\sup_{0\le u\le2}\|T(u)\|$, and let $\delta>0$ correspond to $\eps=1$ in Definition \ref{LMdef}. For $0\le u\le 2$, since $z_n\in D(A)$,
		the map $r\mapsto T(r)z_n$ is continuously differentiable, with
		\[
		\frac{d}{dr}T(r)z_n=T(r)Az_n.
		\] It follows that
		\begin{equation}\label{Tuzneq}
			\|T(u)z_n-z_n\|=\left\|\int_0^uT(r)Az_n\,dr\right\|\le 2K\|Az_n\|.
		\end{equation}
		Fix $n\ge 1$, put $x_n=T(1)z_n$, $y_n=z_n$ and let  $s\in {\cal C}$ with
		$$s(t)=\begin{cases} 2t&\mbox{for}\quad 0\le t\le1,\\
			t+1&\mbox{for}\quad t\ge1.
		\end{cases}$$
			For $t\ge1$ we have $T(t)x_n=T(t+1)z_n=T(s(t))y_n$. For $0\le t\le1$ we have $t+1=2t+(1-t)$, hence
		\begin{equation*}
		T(t)x_n-T(s(t))y_n=T(2t)(T(1-t)z_n-z_n).
		\end{equation*}
		By \eqref{Tuzneq}
		\begin{equation}\label{Ttxn} \|T(t)x_n-T(s(t))y_n\|\le 2K^2\|Az_n\|\quad\mbox{for}\quad 0\le t\le 1.
		\end{equation} 
		Note that 
		\begin{align*}
		\|x_n\|=\|T(1)z_n\|
		&\ge \|z_n\|-\|T(1)z_n-z_n\|\\
		&=1-\|T(1)z_n-z_n\|\\
		&\ge 1-2K\|Az_n\|,
		\end{align*}
		using \eqref{Tuzneq}.
		Since $Az_n\to0$ and $x_n=T(1)z_n\to z_n$ in the sense that
		$\|x_n-z_n\|\to0$, choose $n$ sufficiently large that
		\[
		\|x_n\|>\frac12,
		\quad
		2K^2\|Az_n\|<\frac{\delta}{2}.
		\]
		The preceding estimates then give
		\[
		\|T(t)x_n-T(s(t))y_n\|<\delta\|x_n\|
		\quad\mbox{for}\quad t\ge0.
		\]
		Furthermore, for $t \ge 1$, the semigroup law gives 		\[
		T(t)x_n=T(t)T(1)z_n
		=T(t+1)z_n
		=T(s(t))y_n.
		\]
		It follows that
		\begin{equation*}
			\|T(t)x_n-T(s(t))y_n\|<\delta\|x_n\| \quad\mbox{for all}\quad  t\ge 0,
		\end{equation*}  and positive expansiveness gives $t_*\in[0,1)$ with $y_n=T(t_*)x_n$ or $x_n=T(t_*)y_n$. In the first case $z_n=T(t_*+1)z_n$ with $t_*+1>0$, so $z_n$ is a nonzero periodic point, contradicting Theorem \ref{nethm}\,(ii). In the second case $T(1)z_n=T(t_*)z_n$, i.e.\ $T(t_*)(T(1-t_*)z_n-z_n)=0$; by injectivity (Theorem \ref{nethm}\,(i)) $T(1-t_*)z_n=z_n$ with $1-t_*>0$, again a nonzero periodic point. Both cases are impossible. Thus $0\notin\sigma_{\ap}(A)$.
	\end{proof}

	\subsection{Positive expansiveness and the class of reparametrizations}
		
		In this section, we show that continuous reparametrization in 
		Definition~\ref{LMdef} can be replaced by increasing homeomorphisms. A continuous reparametrization is first replaced by its running maximum, which is nondecreasing, and the latter is then approximated by an increasing homeomorphism; the next lemma carries out this last step. 	
	
	Let  ${\cal H}\subset{\cal C}$ denote the set of increasing homeomorphisms $h:[0,\infty)\to [0,\infty)$ with $h(0)=0$.  		
	\begin{lemma}\label{strictlm} Let $y\in X$ and $\delta>0$ and let $z:[0,\infty)\to[0,\infty)$ be continuous and nondecreasing with $z(0)=0$ and $z(t)\to\infty$ as $t\to\infty$. There exists $h\in{\cal H}$ with 
		\begin{equation}\label{Tht} \|T(h(t))y-T(z(t))y\|<\delta\quad\mbox{for}\quad t\ge 0.
	\end{equation}
	\end{lemma}
	
	\begin{proof}
		Since $u\mapsto T(u)y$ is continuous on $[0,\infty)$, it is uniformly continuous on every compact time interval. Subdividing each interval $[k,k+1]$ into finitely many sufficiently small intervals gives a locally finite sequence
		\[
		0=u_0<u_1<u_2<\cdots,\quad u_j\to\infty,
		\]
		such that
		\begin{equation}\label{Tuv}
		\|T(u)y-T(v)y\|<\delta
		\quad\text{whenever }u,v\in[u_j,u_{j+1}].
		\end{equation}
		For each $j\ge0$, define the first hitting time
		\[
		t_j=\min\{t\ge0:z(t)=u_j\}.
		\]
		These times do exist because $z$ is continuous, $z(0)=0$ and $z(t)\to\infty$ as $t\to\infty$ by the hypothesis. Moreover, since $z$ is nondecreasing, $t_j\le t_{j+1}$, and $t_j=t_{j+1}$ would imply $u_j=u_{j+1}$. Therefore, 
		\[
		0=t_0<t_1<t_2<\cdots.
		\]
		Also $t_j\to\infty$ since otherwise the values $z(t_j)=u_j$ would be unbounded on a compact time interval, contradicting continuity.
		
		Define $h$ by linear interpolation:
		\[
		h(t)=u_j+
		\frac{u_{j+1}-u_j}{t_{j+1}-t_j}(t-t_j),
		\quad t\in[t_j,t_{j+1}].
		\]
		Then $h$ is continuous and strictly increasing, with
		\[
		h(0)=0,\quad h(t_j)=u_j\to\infty.
		\]
		Consequently, it is a homeomorphism of $[0,\infty)$ onto itself.
		
		Finally, if $t\in[t_j,t_{j+1}]$, then monotonicity of $z$ gives
		\[
		u_j=z(t_j)\le z(t)\le z(t_{j+1})=u_{j+1}.
		\]
		By construction, $h(t)\in[u_j,u_{j+1}]$ as well. Applying \eqref{Tuv}, we obtain
		\[
		\|T(h(t))y-T(z(t))y\|<\delta
		\quad\mbox{for}\quad t\ge 0.
		\]
		This proves \eqref{Tht}.
	\end{proof}
	
	\begin{theorem}\label{mthm0}
		Let $T$ be a $C_0$-semigroup on a complex Banach space $X$. The following assertions are equivalent:
		
		(i) $T$ is positively expansive;
		
		(ii) for every $\eps>0$ there exists $\delta>0$ such that, whenever $x,y\in X$ and $s\in {\cal H}$ satisfy
		\begin{equation}\label{Hdn}
			\|T(t)x-T(s(t))y\|<\delta\|x\|\quad\mbox{for}\quad t\ge 0,
		\end{equation}
		there is $t_*\in[0,\eps)$ for which $y=T(t_*)x$ or $x=T(t_*)y$.
	\end{theorem}

	\begin{proof} (i)~$\Rightarrow$~(ii). This follows from $\mathcal H\subset\mathcal C$.
		
		(ii)~$\Rightarrow$~(i). Let $T$ satisfy (ii). Since the identity $s(t)=t$ belongs to ${\cal H}$, condition (ii) implies that $T$ is positively kinematic expansive; hence, by Theorem~\ref{ubothm}, every nonzero orbit of $T$ is unbounded. Moreover, the reparametrization
		$$s(t)=\begin{cases} 2t&\mbox{for}\quad 0\le t\le1,\\ t+1&\mbox{for}\quad t\ge1\end{cases}$$
		used in the proof of Theorem~\ref{nAxthm} belongs to ${\cal H}$, and that proof uses positive expansiveness only for this $s$ and through Theorem~\ref{nethm}. Hence it applies verbatim under (ii) and gives $0\notin\sigma_\ap(A)$. Let $\kappa=\min\{\kappa_1,1\}$, where $\kappa_1$ denotes the constant corresponding to $\eps=1$ in \eqref{cleqsup}:
		\begin{equation}\label{kappa}
			\sup_{0\le t\le 1}\|T(t)x-x\|\ \ge\ \kappa\|x\|\quad\mbox{for}\quad x\in X.
		\end{equation}
		Define 
		\begin{equation}\label{K} 
		K=\sup_{0\le r\le2}\|T(r)\|\ge1,
		\quad C=\frac{4K^2}{\kappa}+3.
		\end{equation}
		Given $\varepsilon>0$, let $\eta>0$ be the constant corresponding to $\varepsilon$ in (ii), and set $\delta=\eta/(C+1)$. We show that this $\delta$ satisfies Definition~\ref{LMdef} for $\varepsilon$.
		
		Suppose $x,y\in X$ and $s\in\mathcal C$ satisfy
		\begin{equation}\label{eq1}
		\|T(t)x-T(s(t))y\|<\delta\|x\|\quad\mbox{for}\quad t\ge0.
		\end{equation}
		The strict inequality in \eqref{eq1} excludes $x=0$.
		
		\textit{Step 1. Construct a nondecreasing reparametrization.}
		
		Note that $s$ is unbounded. Indeed, if $R:=\sup_{t\ge0}s(t)<\infty$, then \eqref{eq1} gives
		\[
		\sup_{t\ge0}\|T(t)x\|\le \delta\|x\|+\sup_{0\le r\le R}\|T(r)\|\,\|y\|<\infty,
		\]
		contradicting the unboundedness of the orbit of $x\ne0$; see Theorem \ref{ubothm}. Define \[z(t)=\max_{0\le u\le t}s(u),\, t\ge 0,\]
		the running maximum of $s$. 		  
		Then, $z:[0,\infty)\to[0,\infty)$ is continuous, nondecreasing, $z(0)=0$ and $z(t)\to\infty$ as $t\to\infty$.

		\textit{Step 2.}\textit{ Prove the claim }
		\begin{equation}\label{claimeq}
			\|T(t)x-T(z(t))y\|<C\delta\|x\|\quad\mbox{for}\quad t\ge0.
		\end{equation}
			
		Let $t\ge 0$. If $s(t)=z(t)$, then
		\[\|T(t)x-T(z(t))y\|=\|T(t)x-T(s(t))y\|<\delta\|x\|<C\delta\|x\|. \]		
		Now, consider $s(t)<z(t)$. Let $(\tau,\tau')$ be the maximal interval containing $t$ such that $s(r)<z(r)$ for all $r\in (\tau, \tau')$. The running maximum $z$ is constant on $(\tau,\tau')$.
		Indeed, if $\tau<a<b<\tau'$ and $z(a)<z(b)$, the maximum
		of $s$ on $[0,b]$ is attained at some $c\in(a,b]$.
		Then $s(c)=z(b)=z(c)$, contradicting $s<z$ on
		$(\tau,\tau')$. So, $z(r)=z(t)=:m$ for all $r\in (\tau,\tau')$. 
		Since $z$ is unbounded, $\tau'<\infty$. By continuity of $z$, $z(\tau)=z(\tau')=m$; by maximality of $(\tau,\tau')$ and continuity of $s-z$, $\tau$ and $\tau'$ do not belong to the open set $\{s<z\}$, so
		$s(\tau)=z(\tau)=z(\tau')=s(\tau')=m$. By \eqref{eq1}
		\begin{equation}\label{Ttaux}
			\|T(\tau)x-T(m)y\|<\delta\|x\|.
		\end{equation}
		 Set $w=\tau'-\tau$.
		For $0\le \nu\le\min\{1,w\}$, set $q=\min\{\nu,\tau\}$, and define
		\[
		d(p)=s(\tau-p)-s(\tau+\nu-p),
		\quad 0\le p\le q.
		\]
		 Recall that $s(u)\le m$ for $0\le u\le\tau'$. Since $\tau+\nu\in[\tau,\tau']$,
		\[
		d(0)=s(\tau)-s(\tau+\nu)\ge0.
		\]
		If $\nu\le\tau$, then
		\[
		d(q)=d(\nu)=s(\tau-\nu)-s(\tau)\le0.
		\]
		If $\nu>\tau$, then
		\[
		d(q)=d(\tau)=s(0)-s(\nu)\le0.
		\]
		Therefore continuity gives $p\in[0,q]$ with $d(p)=0$, or equivalently $s(\tau-p)=s(\tau+\nu-p)$.	
		Applying \eqref{eq1} gives 
		\[
		\|T(\tau+\nu-p)x-T(s(\tau+\nu-p))y\|<\delta\|x\|,
		\]
		and
		\[
		\|T(\tau-p)x-T(s(\tau-p))y\|<\delta\|x\|.
		\]
		It follows from 
		 $s(\tau-p)=s(\tau+\nu-p)$ that
		\[
		\begin{aligned}
			\|T(\tau+\nu-p)x-T(\tau-p)x\|<2\delta\|x\|.
		\end{aligned}
		\]
		Applying $T(p)$ yields
		\begin{align}\label{Ttaug}
			\|T(\tau+\nu)x-T(\tau)x\| \notag 
			&=\|T(p)(T(\tau+\nu-p)x-T(\tau-p)x )\|\\ \notag
			&\le \|T(p)\|\|T(\tau+\nu-p)x-T(\tau-p)x\|\\
			&<2K\delta\|x\|;
		\end{align}
		recall $K$ from \eqref{K}. 
		
		If $w<1$, use $\nu=t-\tau$ in \eqref{Ttaug}, together with \eqref{Ttaux} to obtain
		\[
		\|T(t)x-T(m)y\|<(2K+1)\delta\|x\|<C\delta\|x\|;
		\]
		recall $C$ from \eqref{K}. It remains to consider $w\ge1$. Rewrite \eqref{Ttaug}:
		\[\|T(\nu)T(\tau)x-T(\tau)x\|<2K\delta\|x\|.\]
		Apply \eqref{kappa} to get
		\[\kappa\|T(\tau)x\| \le \sup_{\nu\in[0,1]}\|T(\nu)T(\tau)x-T(\tau)x\|\le 2K\delta\|x\|,\]
		equivalently,
		\begin{equation}\label{Taux} 
			\|T(\tau)x\|\le\frac{2K}{\kappa}\delta\|x\|.
		\end{equation}
		It follows from \eqref{Ttaux} that
		\begin{equation}\label{Ttaux2}
			\|T(m)y\|<\left(\frac{2K}{\kappa}+1\right)\delta\|x\|.
		\end{equation}
		Let $v\in[\tau,t]$ be the first time with $s(v)=s(t)$.
		We show that
		\begin{equation}\label{Tmy}
			\|T(v)x\|\le\frac{2K^2}{\kappa}\delta\|x\|.
		\end{equation}
		Since $s(\tau)=m>s(t)$, continuity and the choice of $v$ imply
		$s(r)>s(v)$ for every $r\in[\tau,v)$.
		If $v-\tau\le1$, then by the definition of $K$ and  \eqref{Taux} 
		\[\|T(v)x\|=\|T(v-\tau)T(\tau)x\|\le \|T(v-\tau)\| \|T(\tau)x\|\le\frac{2K^2}{\kappa}\delta\|x\|.\] 
		If $v-\tau>1$, then $u:=v-1\in(\tau,v)$, and hence
		$s(u)>s(v)$.  Consider 
		\[
		f(\sigma)=\sigma-u+s(u)-s(\sigma),
		\quad \sigma\in[u,v].
		\]
		We have $f(u)=0$ and $f(v)=1+s(u)-s(v)>1.$ Therefore, for each $r\in[0,1]$, there exists $\sigma\in[u,v]$ such that $f(\sigma)=r$, or equivalently,
		\[
		v+r-\sigma+s(\sigma)=1+s(u).
		\]
		Using the semigroup property, $v=u+1$, and this equality, we obtain
		\[
		\begin{aligned}
			&T(v+r-\sigma)\bigl(T(\sigma)x-T(s(\sigma))y\bigr)
			-T(1)\bigl(T(u)x-T(s(u))y\bigr)\\
			&\quad=
			T(v+r)x-T\bigl(v+r-\sigma+s(\sigma)\bigr)y
			-T(u+1)x+T\bigl(1+s(u)\bigr)y\\
			&\quad=T(v+r)x-T(v)x\\
			&\quad=(T(r)-I)T(v)x.
		\end{aligned}
		\]
		Moreover, $
		0\le v+r-\sigma\le 2.$		
		Consequently, \eqref{eq1} and the definition of $K$ imply
		\[
		\|(T(r)-I)T(v)x\|<2K\delta\|x\|,
		\quad 0\le r\le 1.
		\]
		Applying \eqref{kappa} to $T(v)x$, we obtain
		\[
		\kappa\|T(v)x\|
		\le
		\sup_{0\le r\le1}\|(T(r)-I)T(v)x\|
		\le2K\delta\|x\|.
		\]
		Thus
		\[
		\|T(v)x\|\le\frac{2K}{\kappa}\delta\|x\|
		\le\frac{2K^2}{\kappa}\delta\|x\|.
		\]
		This proves  \eqref{Tmy} in both cases.
				
		Since $s(v)=s(t)$, we have
		\[ \|T(t)x-T(v)x\|\le \|T(t)x-T(s(t))y\|+\|T(v)x-T(s(v))y\|<2\delta\|x\|\]
by \eqref{eq1}. Combining this with \eqref{Ttaux2} and \eqref{Tmy},
		\[
			\|T(t)x-T(z(t))y\|
			<\left(\frac{2K^2+2K}{\kappa}+3\right)\delta\|x\|\
			\le C\delta\|x\|,
		\]
		which means that \eqref{claimeq} holds.
		
		\textit{Step 3. Replace the nondecreasing reparametrization by an increasing homeomorphism.}
		
		Using Lemma~\ref{strictlm}, there exists  $h\in\mathcal H$ such that
		\[ \|T(h(t))y-T(z(t))y\|<\delta\|x\|\quad\mbox{for}\quad t\ge0.\]
		Therefore, it follows from \eqref{claimeq} that
		\begin{align*}
		\|T(t)x-T(h(t))y\|
		&\le \|T(t)x-T(z(t))y\|+\|T(z(t))y-T(h(t))y\|\\
		&<(C+1)\delta\|x\|=\eta\|x\|\quad\mbox{for}\quad t\ge0.
		\end{align*}
	Since $h\in{\cal H}$, condition (ii), applied with this $\eta$, supplies $r\in[0,\varepsilon)$ such that $y=T(r)x$ or $x=T(r)y$. As $x,y\in X$ and $s\in{\cal C}$ satisfying \eqref{eq1} were arbitrary, $T$ is positively expansive, which is (i).
	\end{proof}

	\section{Semigroups with uniform exponential expansion}\label{sec4} 
		
	This section establishes uniform exponential expansion as a sufficient condition for positive expansiveness without assuming surjectivity. We then characterize this growth property using the minimum modulus, the lower growth bound, and the approximate point spectrum of a semigroup operator. We then relate it to exponential decay of the inverse semigroup in the invertible case, which recovers the sufficient condition of Lee and Morales, and conclude with a Lyapunov-type characterization on Hilbert spaces and a perturbation result.
	\subsection{Sufficiency and spectral characterization}
	
	\begin{definition}\rm
		Let $T$ be a $C_0$-semigroup on a Banach space $X$. We say that
		$T$ is \textit{uniformly exponentially expanding} if 
		there exist $M>0,a>0$ such that
		\begin{equation}\label{mthmeq}
			\|T(t)x\|\ge Me^{at}\|x\|\quad\mbox{for}\quad x\in X,\ t\ge 0.
		\end{equation}
We also call \eqref{mthmeq} a \textit{uniform exponential expansion}. 
	\end{definition}
	The main result in this section is the following.
	\begin{theorem}\label{mainthm}
		Let $T$ be a $C_0$-semigroup on a Banach space $X$. If $T$ is uniformly exponentially expanding, then it is positively expansive.
	\end{theorem}
	
	\begin{proof} Let $T$ be uniformly exponentially expanding. 
		We claim that for every $\eps>0$ there is $\eta>0$ such that
		\begin{equation}\label{cleq}
			\|T(r)z-z\|\ge\eta\|z\|\quad\mbox{for}\quad z\in X,\ r\ge\eps.
		\end{equation}
		In other words, after a time of at least $\varepsilon$, a nonzero vector cannot return arbitrarily close to itself relative to its norm.
		Fix $\eps>0$ and choose $R>\eps$ with $Me^{aR}>1$. For $r\ge R$, the reverse triangle inequality and \eqref{mthmeq} give
		\begin{equation}\label{cl1}
		\|T(r)z-z\|\ge\|T(r)z\|-\|z\|\ge(Me^{ar}-1)\|z\|\ge(Me^{aR}-1)\|z\| \quad\mbox{for}\quad z\in X.
		\end{equation}
		Now let $\eps\le r\le R$. Choose an integer $n\ge1$ with $Me^{an\eps}>1$ and put $K=\displaystyle\sup_{0\le u\le nR}\|T(u)\|<\infty$. The identity
		\[ 	T(nr)z-z =\sum_{k=0}^{n-1}T(kr)\bigl(T(r)z-z\bigr)\]
		gives
		\[
		\|T(nr)z-z\|\le nK\|T(r)z-z\|.
		\] Moreover,
		\begin{equation*}
		\|T(nr)z-z\|\ge(Me^{anr}-1)\|z\|\ge(Me^{an\eps}-1)\|z\|,
		\end{equation*}
		using \eqref{mthmeq}. It follows that
		\begin{equation}\label{cl2} 
			\|T(r)z-z\|\ge \frac{Me^{an\eps}-1}{nK}\|z\|\quad \mbox{for}\quad z\in X,\ r\in [\eps,R]. 
		\end{equation} 
		Combining \eqref{cl1} and \eqref{cl2}, inequality \eqref{cleq} holds with
		\begin{equation*}
		\eta=\min\left\{Me^{aR}-1,\ \frac{Me^{an\eps}-1}{nK}\right\}>0 .
		\end{equation*}
		
		Fix $\eps>0$, take $\eta$ satisfying \eqref{cleq} and fix $0<\delta<\min\{\tfrac12,\tfrac\eta2\}$. Let $x,y\in X$ and $s\in {\cal C}$ satisfy 
		\begin{equation}\label{LMdefeq1}
			\|T(t)x-T(s(t))y\|<\delta\|x\|\quad\mbox{for}\quad t\ge 0.
		\end{equation} The strict inequality excludes $x=0$, and dividing $x$ and $y$ by $\|x\|$ we may assume $\|x\|=1$. Putting $t=0$ in \eqref{LMdefeq1},
		\begin{equation}\label{x-yeq}
			\|x-y\|<\delta,\quad \|y\|>1-\delta>\tfrac12 .
		\end{equation}
		
		Note that $s(t)\to\infty$ as $t\to\infty$, since otherwise there are $R_0<\infty$ and $\{t_j\}_{j\ge 1}\subset (0,\infty)$ with $t_j\to\infty$ as $j\to\infty$ such that $s(t_j)\le R_0$, and
		\begin{equation*}
		\|T(t_j)x\|\le\|T(s(t_j))y\|+\delta\le\displaystyle\sup_{0\le u\le R_0}\|T(u)\|\|y\|+\delta,
		\end{equation*}
		contradicting $\|T(t_j)x\|\ge Me^{at_j}\to\infty$ as $j\to\infty$. 
		
		Set $d(t)=s(t)-t,\,t\in [0,\infty)$. If $d(t)\ge 0$, then $s(t)=t+d(t)$ and $T(s(t))y=T(t)T(d(t))y$. By linearity, \[T(t)x-T(s(t))y=T(t)(x-T(d(t))y).\] Applying \eqref{mthmeq} to the vector $x-T(d(t))y$ and using \eqref{LMdefeq1} with $\|x\|=1$, one has
		\begin{equation*}
		Me^{at}\|x-T(d(t))y\|\le \|T(t)x-T(s(t))y\|<\delta,
		\end{equation*}
		that is,
		\begin{equation}\label{x-Teq}
			\|x-T(d(t))y\|<\delta M^{-1}e^{-at}.
		\end{equation}
		If $d(t)\le0$, then $t=s(t)+(-d(t))$ with $-d(t)\ge 0$, so $T(t)x=T(s(t))T(-d(t))x$ and \[ T(t)x-T(s(t))y=T(s(t))(T(-d(t))x-y).\]
		 Exactly as before, \[Me^{as(t)}\|T(-d(t))x-y\|<\delta,\] that is,
		\begin{equation}\label{Tdxyeq}
			\|T(-d(t))x-y\|<\delta M^{-1}e^{-as(t)}.
		\end{equation}
		Note that the function $d$ is bounded on $[0,\infty)$. Indeed, otherwise, there is a sequence $t_j\to\infty$ such that either $d(t_j)\to\infty$ or $d(t_j)\to -\infty$. If $d(t_j)\to\infty$, then \eqref{x-Teq} gives \[\|T(d(t_j))y\|\le1+\delta/M,\]  whereas 
		\[\|T(d(t_j))y\|\ge Me^{ad(t_j)}\|y\|\to\infty\]
		as $j\to\infty$. 
		 If $d(t_j)\to-\infty$, then \eqref{Tdxyeq} gives \[\|T(-d(t_j))x\|\le\|y\|+\delta/M,\]
		  whereas \[\|T(-d(t_j))x\|\ge Me^{-ad(t_j)}\to\infty\]
		  as $j\to\infty$. Both cases are impossible.
		
		Since the function $d$ is bounded on $[0,\infty)$, there is a sequence $t_j\to\infty$ and  $d(t_j)\to h$ for some $h\in\R$. By passing to subsequences if necessary, we may assume that 
		either \[ d(t_j) \ge 0\quad\mbox{for all}\quad j=1,2\dots\]
		or  \[ d(t_j) \le 0\quad\mbox{for all}\quad j=1,2\dots\]
		In the first case $h\ge 0$. The orbit $u\mapsto T(u)y$ is continuous, so $T(d(t_j))y\to T(h)y$ as $j\to\infty$, while the right-hand side of \eqref{x-Teq} tends to $0$ as $j\to\infty$. Hence $\|x-T(h)y\|=0$, i.e.\ $x=T(h)y$. In the second case $h\le0$, similarly $T(-d(t_j))x\to T(-h)x$ as $j\to\infty$. The right-hand side of \eqref{Tdxyeq} tends to $0$ because $s(t_j)\to\infty$ (as proved above). Hence $y=T(-h)x$. With $r:=|h|$ we have obtained
		\begin{equation*}
		y=T(r)x\quad\text{or}\quad x=T(r)y .
		\end{equation*}
		It remains to show that $r\in [0,\eps)$. Suppose on the contrary that $r\ge\eps$. If $y=T(r)x$, then \eqref{cleq} gives $\|x-y\|\ge\eta\|x\|=\eta$, contradicting \eqref{x-yeq} since $\delta<\eta$. If $x=T(r)y$, then \eqref{cleq} and \eqref{x-yeq} give $\|x-y\|\ge\eta\|y\|>\eta/2>\delta$, again a contradiction. Hence $r\in[0,\eps)$, which proves positive expansiveness.
	\end{proof}
	\begin{remark}\rm  The lower bound  \[\|T(t)\| \ge Me^{at}\quad\mbox{for}\quad t\ge 0\]
		for some $a,M>0$
		does not imply positive expansiveness. Indeed, let $X=\C^2$ with the Euclidean norm, 
		define
		\[
		T(t)=
		\begin{pmatrix}
			e^{a t}&0\\
			0&1
		\end{pmatrix},
		\quad t\ge0,
		\]
		where $a>0$. Then
		\[
		\|T(t)\|=e^{at}\quad\mbox{for}\quad t\ge 0.
		\]
		However, for $x=(0,1)\in \C^2$, $T(t)(x)=x$ for all $t\ge 0$. This means $T$ has nonzero fixed points, so it is not positively expansive by Theorem \ref{nethm}\,(ii).   
	\end{remark}
	
	Condition \eqref{mthmeq} involves all times and all vectors. The next result shows that it can be checked on a single operator of the semigroup, either through its minimum modulus or through its approximate point spectrum. It is a one sided, continuous time analogue of Hedlund's characterization of uniformly expansive automorphisms in \cite{Hed}.
	
	\begin{theorem}\label{uepeqthm}
		For a $C_0$-semigroup $T$ the following assertions are equivalent:
		\begin{enumerate}
			\item[(i)] $T$ is uniformly exponentially expanding; %i.e.\ \eqref{mthmeq} holds for some $M,a>0$;
			\item[(ii)] $\tn{T(\tau)}>1$ for some $\tau>0$;
			\item[(iii)] $\omega_-(T)>0$;
			\item[(iv)] $r_-(T(\tau))>1$ for some, equivalently every, $\tau>0$;
			\item[(v)] $\sigma_\ap(T(\tau))\subset\{\lambda\in\C:|\lambda|>1\}$ for some, equivalently every, $\tau>0$.
		\end{enumerate}
		Each of these conditions implies that $T$ is positively expansive.
	\end{theorem}
	
	\begin{proof}
		(i)~$\Rightarrow$~(ii). Assume (i) and choose $M,a>0$ as in  \eqref{mthmeq}. Taking the minimum over all unit vectors gives $\tn{T(t)}\ge Me^{at}$. For $\tau>0$ sufficiently large, $Me^{a\tau}>1$, so (ii) follows.
		
		(ii)~$\Rightarrow$~(i). Suppose (ii) holds, and set \[q:=\tn{T(\tau)}>1,\quad K:=\displaystyle\sup_{0\le u\le\tau}\|T(u)\|<\infty.\] For $0\le r\le\tau$ and $x\in X$,
		\begin{equation*}
		q\|x\|\le\|T(\tau)x\|=\|T(\tau-r)T(r)x\|\le K\|T(r)x\|,
		\end{equation*}
		so $\tn{T(r)}\ge q/K$ on $[0,\tau]$. For $t=n\tau+r$ with $n\in\N\cup\{0\}$ and $0\le r<\tau$, the first inequality in \eqref{L1L2} gives \[\tn{T(t)}\ge\tn{T(\tau)}^n\tn{T(r)}\ge q^{n+1}/K.\]
		 Put $a:=(\log q)/\tau>0$. Since $t<(n+1)\tau$, we have $q^{n+1}=e^{a(n+1)\tau}\ge e^{at}$, and therefore
		\begin{equation*}
			\tn{T(t)}\ \ge\ \frac1K\exp\left(\frac{\log q}{\tau}t\right)\quad\mbox{for}\quad t\ge 0,
		\end{equation*}
		which is \eqref{mthmeq} with $M=K^{-1}$ and $a=(\log q)/\tau$, and hence (i) holds.
		
		(ii)~$\Leftrightarrow$~(iii). By Lemma \ref{lgblm}\,(ii), $\omega_-(T)=\displaystyle\sup_{t>0}\frac1t\log\tn{T(t)}$, which is positive if and only if $\log\tn{T(\tau)}>0$ for some $\tau>0$. Equivalently, $\tn{T(\tau)}>1$ for some $\tau>0$. 
		
		(iii)~$\Leftrightarrow$~(iv). Using Lemma \ref{lgblm}\,(iii), $r_-(T(\tau))=e^{\tau\omega_-(T)}$ for every $\tau>0$. Therefore, $r_-(T(\tau))>1$ if and only if $\omega_-(T)>0$.
		
		(iv)~$\Leftrightarrow$~(v). Due to \eqref{MZeq}, $r_-(T(\tau))=\min\{|\lambda|:\lambda\in\sigma_\ap(T(\tau))\}$, the minimum being taken over a nonempty compact set. Hence $r_-(T(\tau))>1$ if and only if all points of $\sigma_\ap(T(\tau))$ have modulus greater than $1$.
		
		Finally, (i) implies positive expansiveness of the semigroup by Theorem \ref{mainthm}.
	\end{proof}
	
	\begin{corollary}\label{isoco}
		Let $V$ be a $C_0$-semigroup of linear isometries on $X$, i.e. $\|V(t)x\|=\|x\|,
		\ x\in X,\ t\ge0$. For every $a>0$ the semigroup $T(t)=e^{at}V(t)$ is positively expansive.
	\end{corollary}
	
	\begin{proof} Since $V(0)=I$, we have $T(0)=I$. For $t,s\ge0$,
		\begin{equation*}
			T(t)T(s)=e^{at}e^{as}V(t)V(s)
			=e^{a(t+s)}V(t+s)=T(t+s).
		\end{equation*}
		Strong continuity of $T$ follows from continuity of the scalar factor and strong
		continuity of $V$. Moreover, the isometry property gives
		\begin{equation*}
			\|T(t)x\|=e^{at}\|V(t)x\|
			=e^{at}\|x\|\quad\mbox{for}\quad x\in X,\ t\ge0.
		\end{equation*}
		Thus the hypothesis of Theorem~\ref{mainthm} holds with
		$M=1$, and $T$ is positively expansive.
	\end{proof}
	
	\begin{example}[Weighted right translations]\label{wrte} \rm 
		Let $X=L^p(0,\infty)$ with $1\le p<\infty$ and let
		\begin{equation*}
		(V(t)f)(r)=\begin{cases}0&\mbox{for}\quad 0<r<t,\\ f(r-t)&\mbox{for}\quad r\ge t.\end{cases}
		\end{equation*}
		It is easy to check that $V$ is a semigroup, which is called the right-translation semigroup on $X$. We have \[\|V(t)f\|_p^p=\int_t^\infty|f(r-t)|^p\,dr=\|f\|_p^p.\] If $\tilde f$ denotes the extension of $f$ by zero to $\R$, then $V(t)f$ is the restriction of $\tilde f(\cdot-t)$ to $(0,\infty)$, so continuity of translations in $L^p(\R)$ gives $\|V(t)f-f\|_p\to 0$ as $t\to 0^+$. Thus $V$ is a $C_0$-semigroup of isometries, and for $a>0$ the semigroup $T(t)=e^{at}V(t)$ is positively expansive by Corollary \ref{isoco}. Every function in $\ran T(t)$ vanishes almost everywhere on $(0,t)$, so $T(t)$ is not surjective for any $t>0$. The inverse semigroup therefore does not exist, and criteria based on exponential decay of the inverse semigroup do not apply.
	\end{example}
	
	\begin{remark}\rm \label{lasrota} In \cite[Section 4]{LM}, Lee and Morales consider the Lasota equation
		$$\partial_tu+x\,\partial_xu=\lambda u,\quad t\ge 0,\ x\ge 0$$  with $u(0,x)=v(x)$ for some $\lambda>0$. Its solution generates the semigroup $T(t)v(x)=e^{\lambda t}v(xe^{-t})$. It was proved in \cite[Theorem 4.1]{LM} that this semigroup is positively expansive on ${\cal C}_b[0,\infty)$, the space of bounded continuous functions on $[0,\infty)$, endowed with the supremum norm. The semigroup $T$ is invertible with $T(t)^{-1}v(x)=e^{-\lambda t}v(xe^{t})$ and $\|T(t)^{-1}\|\le e^{-\lambda t}, \ t\ge 0$.
		However, the semigroup $V$ with $V(t)v(x)=v(xe^{-t})$ is not strongly continuous on ${\cal C}_b[0,\infty)$.  Indeed, let $f(x)=\cos x$ and set
		\[
		t_n=\log\left(1+\frac{1}{2n}\right),
		\quad
		\xi_n=(2n+1)\pi,
		\quad n\ge1.
		\]
		Then $t_n\to0^+$ and
		\[
		e^{-t_n}\xi_n=2n\pi.
		\]
		Consequently,
		\[
		\begin{aligned}
			\|V(t_n)f-f\|_\infty
			&\ge |f(e^{-t_n}\xi_n)-f(\xi_n)|\\
			&=|\cos(2n\pi)-\cos((2n+1)\pi)|\\
			&=2.
		\end{aligned}
		\]
		Thus $V$ is not strongly continuous on
		$\mathcal C_b([0,\infty))$, so results on $C_0$-semigroups do not apply to it on this space.
		
		 To obtain a $C_0$-semigroup to which the preceding results apply,  we work on the closed subspace 
		 \[  X={\cal C}_0([0,\infty)) 
		 =\left\{f\in C([0,\infty);\C):
		 \lim_{x\to\infty}f(x)=0\right\}, \]
		 endowed with the supremum norm. 
		 Fix $f\in X$ and $\eps>0$. Choose $R>0$ such that $
		 |f(x)|<\frac{\varepsilon}{4}, \ x\ge R$.
		 For $0\le t\le\log 2$ and $x\ge 2R$, we have $xe^{-t}\ge R$. Therefore,
		 \begin{equation}\label{vd1}
		 |f(xe^{-t})-f(x)|
		 \le |f(xe^{-t})|+|f(x)|
		 <\frac{\varepsilon}{2}.
		 \end{equation}
		 On the other hand, $f$ is uniformly continuous on $[0,2R]$. Choose $\eta>0$ such that
		 whenever  $ u,v\in[0,2R], |u-v|<\eta$
		 \[
		 |f(u)-f(v)|<\frac{\eps}{2}.
		 \]
		 For $x\in [0, 2R]$, we have $xe^{-t}\in [0,2R]$ and $
		 |xe^{-t}-x|\le2R(1-e^{-t})$.  Choose $t_\eps>0$ such that
		 $ t_\varepsilon\le\log2$ and $2R(1-e^{-t_\eps})<\eta$ to have 
		 \begin{equation}\label{vd2}|f(xe^{-t})-f(x)|<\frac{\eps}{2}
		 \end{equation} 
		 for $0\le t<t_\eps$. 
		 		 Combining \eqref{vd1} and \eqref{vd2}, we obtain
		 		 \[|f(xe^{-t})-f(x)|<\frac{\eps}{2}\quad\mbox{for}\quad x\in [0,\infty),\ 0\le t<t_\eps.\]
		This yields
		 \[
		 \|V(t)f-f\|_\infty\le\frac{\varepsilon}{2}<\varepsilon
		 \quad\mbox{for}\quad 0\le t<t_\varepsilon.
		 \]
		 Thus $V(t)f\to f$ in the supremum norm as $t\to 0^+$. Together with the semigroup law and the isometry property, this gives continuity of every orbit at every nonnegative time. Consequently, $V$ is a $C_0$-semigroup.
		 Moreover $x\mapsto xe^{-t}$ is a bijection of $[0,\infty)$, so
		\begin{equation*}
		\|T(t)v\|_\infty=e^{\lambda t}\|v\|_\infty ,
		\end{equation*}
		i.e.\ $T=e^{\lambda t}V$ with $V$ a $C_0$-semigroup of isometries. Corollary \ref{isoco} applies and gives positive expansiveness, on ${\cal C}_0([0,\infty))$. The failure of strong continuity on ${\cal C}_b[0,\infty)$ does not by itself disprove expansiveness there. 
	\end{remark}
	
	\subsection{The invertible case} 
			When the semigroup operators are invertible, lower bounds for the semigroup become upper bounds for its inverse. In this case uniform exponential expansion is equivalent to exponential decay of the inverse semigroup, and the sufficient condition of Lee and Morales \cite{LM} turns out to be a special case.

	Let $T=(T(t))_{t\geq 0}$ be invertible with inverse semigroup $S(t)=T(t)^{-1}$. By  \eqref{L-1}
	\begin{equation}\label{triple}\tn{T(t)}=\|S(t)\|^{-1}.
	\end{equation} 
	For $L\in\LL(X)$, we denote the \textit{spectral radius} by \[r(L)=\max\{|\lambda|:\lambda\in\sigma(L)\},\] which satisfies Gelfand's formula $r(L)=\displaystyle\lim_{n\to\infty}\|L^n\|^{1/n}$, and
	\begin{equation}\label{o0}
		\omega_0(S)=\inf\{\omega\in\R:\ \exists\,C\ge1,\ \|S(t)\|\le Ce^{\omega t}\ \ \forall t\ge 0\}\in[-\infty,\infty)
	\end{equation}
	for the \textit{growth bound}.
	
		The next result lists equivalent forms of uniform exponential expansion for invertible semigroups. A single time suffices: it is enough that one operator of the semigroup has an inverse of norm less than one.
	\begin{proposition}\label{T-1pro}
		For a $C_0$-semigroup $T$ the following assertions are equivalent:
		\begin{enumerate}
			\item[(i)] $T(\tau)$ is invertible and $\|T(\tau)^{-1}\|<1$ for some $\tau>0$;
			\item[(ii)] every $T(t)$ is invertible and $\|T(t)^{-1}\|\to 0$ as $t\to\infty$;
			\item[(iii)] every $T(t)$ is invertible and $e^{bt}\|T(t)^{-1}\|\to 0$ as $t\to\infty$ for some $b>0$;
			\item[(iv)] every $T(t)$ is invertible and the inverse semigroup $S$ satisfies $\omega_0(S)<0$;
			\item[(v)] $\sigma(T(\tau))\subset\{\lambda\in\C:|\lambda|>1\}$ for some $\tau>0$.
		\end{enumerate}
		Each of these conditions implies uniform exponential expansion, hence positive expansiveness. Conversely, if every $T(t)$ is invertible, then uniform exponential expansion implies (i)--(v).
	\end{proposition}
	
	\begin{proof}
		In (i)--(iv), invertibility at one time implies invertibility at all times, and $S$ is a $C_0$-semigroup; see \cite[Theorem I.6.5]{Pazy}.
		
		(i) $\Rightarrow$ (iii). By \eqref{triple},  $\tn{T(\tau)}=\|T(\tau)^{-1}\|^{-1}>1$. It follows from Theorem \ref{uepeqthm} that there exist $M,a>0$ such that \[\tn{T(t)}\ge Me^{at}\quad\mbox{for all}\quad t\ge 0,\] that is, \[\|T(t)^{-1}\|\le M^{-1}e^{-at}\quad\mbox{for all}\quad t\ge 0.\]
		For $b\in(0,a)$ we get \[\lim_{t\to\infty}e^{bt}\|T(t)^{-1}\|\le\lim_{t\to\infty} M^{-1}e^{-(a-b)t}= 0.\]
		
		(iii) $\Rightarrow$ (iv). Using (iii) and local boundedness, there exists $C>0$ such that \[\|S(t)\|\le Ce^{-bt}\quad\mbox{for all}\quad t\ge 0,\] so $\omega_0(S)\le-b<0$, by the definition of $\omega_0$ in \eqref{o0}.
		
		(iv) $\Rightarrow$ (ii). If $\omega_0(S)<\omega<0$, then $\|T(t)^{-1}\|=\|{S(t)}\|\le Ce^{\omega t}\to 0$ as $t\to\infty$.
		
		(ii) $\Rightarrow$ (i). The norm $\|T(t)^{-1}\|<1$  for large $t$.
		
		(i) $\Rightarrow$ (v). For an invertible operator $L$ and $\lambda\ne0$ we have 
		\[\lambda I-L=-\lambda L(\lambda^{-1}I-L^{-1}),\]
		so $\lambda\in\sigma(L)$ if and only if $\lambda^{-1}\in\sigma(L^{-1})$. Thus, if $\|T(\tau)^{-1}\|<1$ for some $\tau$ and $\lambda\in\sigma(T(\tau))$ we have $\lambda^{-1}\in\sigma(S(\tau))$, hence 
		\[|\lambda|^{-1}\le r(S(\tau))\le\|S(\tau)\|<1\]
		and hence $|\lambda|>1$.
		
		(v) $\Rightarrow$ (i). Since $\sigma(T(\tau))\subset\{\lambda\in\C:|\lambda|>1\}$, we have $0\notin\sigma(T(\tau))$ and hence the operator $T(\tau)$ is invertible. The set $\sigma(T(\tau))$ is nonempty and compact, so $d:=\min\{|\lambda|:\lambda\in\sigma(T(\tau))\}>1$, and $r(S(\tau))=d^{-1}<1$. By Gelfand's formula,
		\[
		\lim_{n\to\infty}\|S(\tau)^n\|^{1/n}
		=r(S(\tau))<1,
		\] so $\|S(\tau)^n\|<1$ for large $n$, and $S(\tau)^n=S(n\tau)=T(n\tau)^{-1}$. Thus (i) holds with $n\tau$ in place of $\tau$.
		
		Assume that (i) holds, $\|T(\tau)^{-1}\|<1$ for some $\tau>0$. Then $\tn{T(\tau)}=\|T(\tau)^{-1}\|^{-1}>1$. This is equivalent to the uniform exponential expansion from Theorem \ref{uepeqthm}, and hence $T$ is positively expansive. Conversely, if every $T(t)$ is invertible and uniform exponential expansion holds, then $\|T(\tau)^{-1}\|=\tn{T(\tau)}^{-1}<1$ by Theorem \ref{uepeqthm}, which is (i) and hence this implies (i)-(v).
	\end{proof}

	\begin{remark}\rm (a)  For invertible semigroups Proposition \ref{T-1pro} and Theorem \ref{uepeqthm} therefore coincide; the gain of Theorem \ref{uepeqthm} lies in semigroups without inverse.
		
		(b)	If the conditions of Proposition \ref{T-1pro} hold with $\sigma(T(\tau))\subset\{|\lambda|>1\}$, the spectral inclusion theorem $e^{\tau\sigma(A)}\subset\sigma(T(\tau))$ (see \cite[Theorem IV.3.6]{EN}) gives $\re\lambda>0$ for every $\lambda\in\sigma(A)$; more precisely,
		\begin{equation*}
			\inf_{\lambda\in\sigma(A)}\re\lambda\ \ge\ \frac1\tau\log\min_{\mu\in\sigma(T(\tau))}|\mu|\ >0 ,
		\end{equation*}
		with convention $\inf\varnothing=+\infty$. This concerns the conditions of the proposition and not positive expansiveness alone; see Remark \ref{diaexp}. {\hfill$\Diamond$}
	\end{remark}
	
	\begin{remark} \rm 
		Lee and Morales \cite[Theorem 1.3(2)]{LM} prove positive expansiveness under the resolvent condition $(-\infty,\alpha)\subset\rho(A)$, $\|(\gamma I-A)^{-1}\|\le(\alpha-\gamma)^{-1}$ for $\gamma<\alpha$. Their route has two steps. First, applying the Hille--Yosida theorem to $\alpha I-A$, they show that $T(t)$ is invertible with $\|T(t)^{-1}\|\le e^{-\alpha t}$ \cite[Lemma 2.5]{LM}. Second, they show that this inverse bound implies positive expansiveness \cite[Lemma 2.6]{LM}. Their hypotheses give $\|T(t)^{-1}\|\le e^{-\alpha t}<1$ for every $t>0$, so they are a special case of Proposition \ref{T-1pro}\,(i). The proposition shows that it suffices to have $\|T(\tau)^{-1}\|<1$ at a single time, and Theorem \ref{uepeqthm} shows that invertibility itself can be dispensed with. For bounded generators, Theorem \ref{Abdeqthm} uses only the location of the spectrum; it yields exponential decay of $T(t)^{-1}$ with a multiplicative constant that need not equal one in the given norm. {\hfill $\Diamond$}
	\end{remark}
	
	\subsection{Lyapunov characterization and robustness}
		
	The following result characterizes uniform exponential expansion on Hilbert spaces by a Lyapunov-type inequality.
	
	\begin{theorem}\label{Hthm1}
		Let $A$ generate a $C_0$-semigroup $T$ on a complex Hilbert space $H$. The following assertions are equivalent:
		\begin{enumerate}
			\item[(i)] $T$ is uniformly exponentially expanding;
			\item[(ii)] there exist a bounded self-adjoint operator $P\in\LL(H)$ and constants $m,L,c>0$ such that
			\begin{equation}\label{mPl}
				m\|x\|^2\le\langle Px,x\rangle\le L\|x\|^2\quad\mbox{for}\quad x\in H
			\end{equation}
			and
			\begin{equation}\label{PAc}
				2\re\langle PAx,x\rangle\ge c\|x\|^2\quad\mbox{for}\quad x\in D(A).
			\end{equation}
		\end{enumerate}
		In either case $T$ is positively expansive.
	\end{theorem}
	
	\begin{proof}
		(ii)~$\Rightarrow$~(i). Fix $x\in D(A)$. Then $T(t)x\in D(A)$ and $\frac{d}{dt}T(t)x=AT(t)x$ (\cite[Lemma II.1.3]{EN}), so $t\mapsto T(t)x$ is continuously differentiable in $H$ and, $P$ being bounded, $E(t)=\langle PT(t)x,T(t)x\rangle$ is continuously differentiable with \[E'(t)=\langle PAT(t)x,T(t)x\rangle+\langle PT(t)x,AT(t)x\rangle.\]
		 As $P$ is self-adjoint, the second term is the complex conjugate of the first. Using \eqref{PAc} with $T(t)x\in D(A)$ and then \eqref{mPl},
		\begin{equation*}
			E'(t)=2\re\langle PAT(t)x,T(t)x\rangle\ \ge\ c\|T(t)x\|^2\ \ge\ \frac cL E(t).
		\end{equation*}
		Hence 
		\[\frac{d}{dt}(e^{-ct/L}E(t))=e^{-ct/L}\left(E'(t)-\frac cLE(t)\right)\ge 0,\]
		 and so $E(t)\ge e^{ct/L}E(0)$. By \eqref{mPl},
		\begin{equation*}
			L\|T(t)x\|^2\ge E(t)\ge e^{ct/L}E(0)\ge me^{ct/L}\|x\|^2 .
		\end{equation*}
		For general $x\in H$ choose $x_j\in D(A)$ with $x_j\to x$. For fixed $t$ we have $T(t)x_j\to T(t)x$ because $T(t)$ is bounded, so the estimate passes to the limit. Taking square roots, we obtain
		\begin{equation*}
			\|T(t)x\|\ge\sqrt{m/L}e^{ct/(2L)}\|x\|\quad\mbox{for}\quad x\in H,\, t\ge 0.
		\end{equation*}
		
		(i)~$\Rightarrow$~(ii). Let \eqref{mthmeq} hold, choose $\tau>0$ with $q:=Me^{a\tau}>1$ and put $K=\displaystyle\sup_{0\le s\le\tau}\|T(s)\|$. The sesquilinear form \[p(x,y)=\int_0^\tau\langle T(s)x,T(s)y\rangle ds\] is well defined because the integrand is continuous in $s$, it is sesquilinear, and by the Cauchy--Schwarz inequality \[|p(x,y)|\le\int_0^\tau\|T(s)x\|\|T(s)y\|ds\le\tau K^2\|x\|\|y\|.\]
		 By the Riesz representation theorem (see \cite[Theorem~II.2.2]{Conway}) for bounded sesquilinear forms there is a unique $P\in\LL(H)$ with $\langle Px,y\rangle=p(x,y)$, and $P$ is self-adjoint because $p(y,x)=\overline{p(x,y)}$. Moreover, by \eqref{mthmeq} and the definition of $K$ we have \[\|T(s)x\|^2\ge M^2e^{2as}\|x\|^2 \quad\mbox{and}\quad \|T(s)x\|^2\le K^2\|x\|^2\quad\mbox{on}\quad [0,\tau],\]
		  so
		\begin{equation*}
			M^2\int_0^\tau e^{2as}ds\|x\|^2\le\langle Px,x\rangle=\int_0^\tau\|T(s)x\|^2ds\le\tau K^2\|x\|^2,
		\end{equation*}
		so \eqref{mPl} holds with $m=M^2(e^{2a\tau}-1)/(2a)$ and $L=\tau K^2$. For $x\in D(A)$, the orbit $u(s)=T(s)x$ is continuously differentiable, with
		\[
		u'(s)=T(s)Ax=AT(s)x.
		\]
		The product rule for the inner product gives
		\[
		\begin{aligned}
			\frac{d}{ds}\|u(s)\|^2
			=\langle u'(s),u(s)\rangle
			+\langle u(s),u'(s)\rangle=2\re\langle u'(s),u(s)\rangle.
		\end{aligned}
		\] Therefore,
		\[
		\frac{d}{ds}\|T(s)x\|^2
		=2\re\langle T(s)Ax,T(s)x\rangle
		\]
		implies \begin{align*}
			2\re\langle PAx,x\rangle&=\int_0^\tau2\re\langle T(s)Ax, T(s)x\rangle ds=\int_0^\tau\frac d{ds}\|T(s)x\|^2ds\\
			&=\|T(\tau)x\|^2-\|x\|^2\ge(q^2-1)\|x\|^2 .
		\end{align*}
		Thus \eqref{PAc} holds with $c=q^2-1>0$.
		
		The final assertion follows from Theorem \ref{mainthm}.
	\end{proof}
	
	The final result in this section shows that uniform exponential expansion is robust under small bounded perturbations of the generator. 
	\begin{proposition}\label{pertprop}
		Let $A$ generate a uniformly exponentially expanding $C_0$-semigroup $T$
		on a Banach space $X$. Then there is $\beta_0>0$ such that for every
		$B\in\LL(X)$ with $\|B\|<\beta_0$, the semigroup $S$ generated by $A+B$
		is uniformly exponentially expanding, hence positively expansive.
	\end{proposition}
	
	\begin{proof}
		By Theorem \ref{uepeqthm}\,(ii) there is $\tau>0$ with
		$q:=\tn{T(\tau)}>1$. Choose $M\ge1$ and $\omega\ge0$ with
		$\|T(t)\|\le Me^{\omega t}$ for $t\ge0$. By the bounded perturbation
		theorem \cite[Theorem~III.1.3]{EN}, $A+B$ generates a $C_0$-semigroup $S$
		with $$\|S(t)\|\le Me^{(\omega+M\|B\|)t}\quad \mbox{for}\quad t\geq 0,$$ and
		\[
		S(t)x=T(t)x+\int_0^t T(t-s)BS(s)x\,ds\quad\mbox{for}\quad x\in X,\ t\ge0.
		\]
		Hence
		\[
		\|S(\tau)-T(\tau)\|
		\le\int_0^\tau Me^{\omega(\tau-s)}\|B\|\,Me^{(\omega+M\|B\|)s}\,ds
		\le \tau M^2\|B\|\,e^{(\omega+M\|B\|)\tau}.
		\]
		For $L_1,L_2\in\LL(X)$ and $\|x\|=1$ we have
		$$\|L_1x\|\ge\|L_2x\|-\|(L_1-L_2)x\|\ge\tn{L_2}-\|L_1-L_2\|,$$
		so $$\tn{L_1}\ge\tn{L_2}-\|L_1-L_2\|.$$ Therefore
		\[
		\tn{S(\tau)}\ \ge\ q-\tau M^2\|B\|\,e^{(\omega+M\|B\|)\tau}.
		\]
		Choose $\beta_0>0$ with $\tau M^2\beta_0e^{(\omega+M\beta_0)\tau}<q-1$.
		If $\|B\|<\beta_0$, then $\tn{S(\tau)}>1$, and
		Theorem \ref{uepeqthm} shows that $S$ is uniformly exponentially
		expanding. Positive expansiveness follows from Theorem \ref{mainthm}.
	\end{proof}
	
	\begin{remark}\rm It is worth mentioning that positive expansiveness itself is not preserved under small bounded
		perturbations. Indeed, let $A$ be the diagonal operator on $\ell^p(\N;\C)$,
		$1\le p<\infty$, with $\lambda_n=\frac1n+i$, which generates a positively
		expansive semigroup by Corollary \ref{diagcor} below. For $n\ge1$ let
		$P_nx=x_ne_n$ and $A_n=A-\frac1nP_n$. Then $\|A_n-A\|=\frac1n\to 0$, but
		$A_ne_n=ie_n$, so $i\in\sigma_p(A_n)$ and $e^{tA_n}$ is not positively
		expansive by Theorem \ref{nethm}\,(iii). Thus uniform exponential expansion
		is an open property under bounded perturbations, while positive
		expansiveness is not.
	\end{remark}
	
	\section{Semigroups with a uniform lower bound}\label{sec5}
	Uniform exponential expansion is sufficient for positive expansiveness, but it is not necessary. In this section we characterize positive expansiveness for semigroups that satisfy only a uniform lower bound. Under this assumption, the two necessary conditions of Section~\ref{sec3}, unboundedness of every nonzero orbit and boundedness from below of the generator, turn out to be sufficient as well. For diagonal semigroups this yields explicit criteria in terms of the diagonal entries, which separate the three notions studied in the paper.
	
	The lower bound prevents orbits from shrinking. Combined with the displacement criterion of Lemma~\ref{eqlm}, it keeps every admissible reparametrization within a short distance of the identity, which is the key step in the following sufficient condition. The main result in this section is the following theorem.
	\begin{theorem}\label{mainthm2}
		Let $T$ be a $C_0$-semigroup with generator $A$ on a Banach space $X$. Suppose that
		\begin{enumerate}
			\item[(a)] there is $m>0$ with $\|T(t)x\|\ge m\|x\|$ for all $x\in X$ and $t\ge 0$;
			\item[(b)] every nonzero vector has an unbounded orbit: $\displaystyle\sup_{t\ge 0}\|T(t)x\|=\infty$ for $x\ne0$;
			\item[(c)] $0\notin\sigma_\ap(A)$.
		\end{enumerate}
		Then $T$ is positively expansive.
	\end{theorem}
	
	\begin{proof} Fix $\eps>0$ and put $\eps_1=\eps/2$. By Lemma \ref{eqlm}, assumption (c) implies that there exists
		$\kappa_1:=\kappa_{\eps_1}$ such that 
		\begin{equation}\label{sup0}
			\sup_{0\le t\le\eps_1}\|T(t)x-x\|\ \ge\ \kappa_1\|x\|\quad\mbox{for}\quad x\in X.
		\end{equation}
		Let $K_1=\displaystyle\sup_{0\le t\le\eps_1}\|T(t)\|$ and put
		$ 
		\delta:=\frac{\kappa_1}{2(K_1+1+1/m)}.
		$ Let $x,y\in X$ and $s\in {\cal C}$ satisfy
		\begin{equation}\label{mthm2eq} 
			\|T(t)x-T(s(t))y\|<\delta\|x\|\quad\mbox{for}\quad t\ge 0.
		\end{equation} 
		As before $x\ne0$ and we may assume $\|x\|=1$; then $\|x-y\|<\delta$ by putting $t=0$. Put $d(t)=s(t)-t,\, t\ge 0$.		
		We claim that $|d(t)|<\eps_1$ for all $t\ge 0$. 
		
		Suppose $d(t_1)\ge\eps_1$ for some $t_1>0$. Since $d$ is continuous and $d(0)=0$, every $r\in[0,\eps_1]$ equals $d(t_r)$ for some $t_r\ge 0$. For such $t_r$ we have $s(t_r)=t_r+r$, and \[T(t_r)x-T(s(t_r))y=T(t_r)(x-T(r)y).\] 
		Then
		\begin{equation*}
			\|x-T(r)y\|\le\frac1m\|T(t_r)x-T(s(t_r))y\|<\frac\delta m, 
		\end{equation*}
		using \eqref{mthm2eq}. This yields
		\begin{equation*}
			\|T(r)x-x\|\le\|T(r)(x-y)\|+\|T(r)y-x\|<K_1\delta+\frac\delta m<\kappa_1\quad\mbox{for}\quad 0\le r\le\eps_1,
		\end{equation*}
		which contradicts \eqref{sup0} applied to the unit vector $x$. 
		
		Suppose instead $d(t_2)\le-\eps_1$ for some $t_2>0$. Then every $\tau\in[0,\eps_1]$ equals $-d(t_\tau)$ for some $t_\tau$, and for such $t_\tau$ we have $t_\tau=s(t_\tau)+\tau$, so \[T(t_\tau)x-T(s(t_\tau))y=T(s(t_\tau))(T(\tau)x-y).\] 
		By (a) and \eqref{mthm2eq}, $\|T(\tau)x-y\|<\delta/m$. Hence
		\begin{equation*}
			\|T(\tau)x-x\|\le\|T(\tau)x-y\|+\|y-x\|<\frac\delta m+\delta<\kappa_1\quad\mbox{for}\quad 0\le \tau\le\eps_1,
		\end{equation*}
		again contradicting \eqref{sup0}. In particular $s(t)\ge t-\eps_1\to\infty$.
		
		Next, we claim that if $u_j\to\infty$, $z_j\to z$ in $X$ and $\|T(u_j)z_j\|<\delta$ for all $j$, then $z=0$. To get this, fixing $u\ge 0$, for $j$ large, (a) applied to $T(u)z_j$ gives
		\begin{equation*}
			\|T(u)z_j\|\le\frac1m\|T(u_j-u)T(u)z_j\|=\frac1m\|T(u_j)z_j\|<\frac\delta m,
		\end{equation*}
		and letting $j\to\infty$ we obtain $\|T(u)z\|\le\frac{\delta}m$. As $u$ is arbitrary, the orbit of $z$ is bounded, and (b) yields $z=0$.
		
		Take $t_j\to\infty$. Since $|d(t_j)|<\eps_1$ for all $j=1,2,\dots,$ a subsequence of $\{d(t_j)\}$ converges to some $h\in[-\eps_1,\eps_1]$, and passing to a further subsequence we may assume that either all $d(t_j)\ge 0$ or all $d(t_j)\le0$. In the first case, by the continuity of the orbit of $y$, 
		\[ z_j=x-T(d(t_j))y\to x-T(h)y\quad\mbox{as}\quad j\to\infty\]
		and
		\[\|T(t_j)z_j\|=\|T(t_j)x-T(s(t_j))y\|<\delta\quad\mbox{for all}\quad j.\]
		Using the claim above, $x=T(h)y$. In the second case,  \[z_j=T(-d(t_j))x-y\to T(-h)x-y\quad\mbox{as}\quad j\to\infty \]
		and \[\|T(s(t_j))z_j\|=\|T(t_j)x-T(s(t_j))y\|<\delta.\] Since  $s(t_j)\to\infty$, the claim above  gives $y=T(-h)x$. In either case one of the orbit relations holds at $t_*=|h|\le\eps_1<\eps$.
	\end{proof}
	
	Conversely, conditions (b) and (c) of Theorem~\ref{mainthm2} are necessary for positive expansiveness, so under (a) they characterize it.
	\begin{theorem}\label{charthm}
		Let $T$ be a $C_0$-semigroup with generator $A$ on Banach space $X$. Suppose that there exists $m>0$ such that 
		\begin{equation}\label{mthm3eq}
			\|T(t)x\|\ge m\|x\|\quad\mbox{for}\quad x\in X,\, t\ge 0.
		\end{equation} 
		Then the following assertions are equivalent:
		\begin{enumerate}
			\item[(i)] $T$ is positively expansive;
			\item[(ii)] every nonzero orbit of $T$ is unbounded, and $0\notin\sigma_\ap(A)$;
			\item[(iii)] $T$ is positively kinematic expansive, and $0\notin\sigma_\ap(A)$.
		\end{enumerate}
	\end{theorem}
	
	\begin{proof}
		(ii)~$\Rightarrow$~(i). This is Theorem \ref{mainthm2}.
		
		(i)~$\Rightarrow$~(ii). Suppose that $T$ is positively expansive. The property $0\notin \sigma_\ap(A)$ is Theorem \ref{nAxthm}. Suppose on the contrary that $x\ne0$ had a bounded orbit. According to Lemma \ref{bolm}, there are $r>0$ and $\alpha\in(0,\infty)\setminus\{1\}$ with $T(r)x=\alpha x$. Then $\|T(nr)x\|=\alpha^n\|x\|$ is bounded only if $\alpha<1$, but then it tends to $0$, contradicting  $\|T(nr)x\|\ge m\|x\|$ by \eqref{mthm3eq}.
		
		(ii)~$\Leftrightarrow$~(iii). This is Theorem \ref{ubothm}.
	\end{proof}
	
	\begin{remark}\rm (a) It is worth mentioning that Theorem \ref{mainthm2} holds on real Banach spaces.
		
		(b)	Theorem \ref{mainthm} is a special case of Theorem \ref{mainthm2}. Indeed, \eqref{mthmeq} implies (a) with $m=M$ (as $e^{at}\ge1$) and (b). To see (c), suppose $Ax_n\to 0$ for unit vectors $x_n\in D(A)$. Due to $\frac{d}{ds}T(s)x_n=T(s)Ax_n$, 		
		\[
		T(t)x_n-x_n=\int_0^tT(s)Ax_n\,ds.
		\]
		For each fixed $t>0$, local boundedness of a $C_0$-semigroup ensures that \[K_t:=\sup_{0\le s\le t}\|T(s)\|<\infty.\]
		Consequently,
		\[
		\|T(t)x_n-x_n\|
		\le tK_t\|Ax_n\|\to 0,
		\]
		whereas \eqref{mthmeq} gives $\|T(t)x_n\|\ge Me^{at}>1$ for large $t$, a contradiction. 
		
		{\hfill$\Diamond$}
	\end{remark}
	
	The following semigroups form a class in which all three notions can be decided explicitly from the diagonal entries. When all entries lie in the closed right half-plane, they satisfy hypothesis (a) of Theorem~\ref{mainthm2} with $m=1$.

	\begin{corollary}[Diagonal semigroups]\label{diagcor}
		Let $X=\ell^p(\N;\C)$ with $1\le p<\infty$, and let $\{\lambda_n\}_{n\ge1}\subset\C$ with $\displaystyle\sup_n\re\lambda_n<\infty$. Let $A$ be the diagonal operator $A\{x_n\}=\{\lambda_nx_n\}$ with its maximal domain, which generates the $C_0$-semigroup $T(t)\{x_n\}=\{e^{\lambda_nt}x_n\}$. Then
		
		(a) $T$ is positively kinematic expansive if and only if \begin{equation}\label{disagcoreqa}\re\lambda_n>0\quad\mbox{for all}\quad n\ge 1;
		\end{equation} 
		
		(b) $T$ is positively expansive if and only if
		\begin{equation}\label{disagcoreq}
			\re\lambda_n>0\quad\text{for all}\quad n\ge 1\quad\text{and}\quad\inf_n|\lambda_n|>0 .
		\end{equation}
	\end{corollary}
	
	\begin{proof}
		Write
		\[
		D(A)=\{x\in X:\{\lambda_nx_n\}_{n\ge1}\in X\},
		\]
		and denote by $e_n$ the standard unit vectors. Then $e_n\in D(A)$, $\|e_n\|=1$, and
		\[
		Ae_n=\lambda_ne_n,\quad T(t)e_n=e^{\lambda_nt}e_n.
		\]
		
		(a) \textit{Positive kinematic expansiveness.} Suppose that $T$ is positively kinematic expansive. If $\re\lambda_n\le0$ for some $n$, then $\|T(t)e_n\|=e^{t\re\lambda_n}\le1$, so the orbit of $e_n$ is bounded, contradicting Theorem \ref{ubothm}. Hence, $\re\lambda_n>0$ for all $n\ge 1$. Conversely, suppose that  $\re\lambda_n>0$ for all $n\ge 1$ and $x=\{x_n\}\ne0$, then $x_n\ne0$ for some $n$ and $$\|T(t)x\|\ge e^{t\re\lambda_n}|x_n|\to\infty \quad\mbox{as}\quad t\to\infty.$$ By Theorem \ref{ubothm}, $T$ is positively kinematic expansive.
		
		(b) \textit{Positive expansiveness.}
		Assume that $T$ is positively expansive. Since $\lambda_n\in \sigma_p(A)$, Theorem \ref{nethm}\,(iii) gives $\re \lambda_n>0$ for $n\ge 1$. 
		Furthermore, positive expansiveness implies $0\notin\sigma_{\mathrm{ap}}(A)$. If $\displaystyle\inf_n|\lambda_n|=0$, there exists a strictly increasing sequence
		$\{n_k\}_{k\ge1}\subset\mathbb N$ such that
		\[
		\|Ae_{n_k}\|=|\lambda_{n_k}|\to 0\quad\mbox{as}\quad k\to\infty.
		\]
		Since $\|e_{n_k}\|=1$ for every $k$, this would imply $0\in\sigma_{\mathrm{ap}}(A)$, a contradiction. Thus $		\inf_n|\lambda_n|>0.$
		
		Conversely, assume that \eqref{disagcoreq} holds. We check the three hypotheses of Theorem \ref{mainthm2}. We first have
		\[
		\|T(t)x\|_p^p
		=\sum_{n=1}^{\infty}e^{pt\re\lambda_n}|x_n|^p
		\ge\sum_{n=1}^{\infty}|x_n|^p
		=\|x\|_p^p.
		\]
		Thus (a) holds with $m=1$. 		
		Next, given $x=\{x_n\}\in X, x\ne0$, choose an index $n$ with $x_n\ne0$. Then 
		\[
		\|T(t)x\|
		\ge e^{t\re\lambda_n}|x_n|
		\to\infty.
		\]
		Hence every nonzero forward orbit is unbounded, so (b) holds.
		Finally, for $x\in D(A)$, 
		\[ 	\|Ax\|_p^p
		=\sum_{n=1}^{\infty}|\lambda_n|^p|x_n|^p
		\ge\kappa^p\|x\|_p^p,
		\]
		where $\kappa=\inf_n|\lambda_n|>0$.
		Therefore, 
		\[
		\|Ax\|\ge\kappa\|x\|,
		\quad x\in D(A),
		\]
		which is equivalent to $0\notin\sigma_{\mathrm{ap}}(A)$.
		
		All three hypotheses of Theorem \ref{mainthm2} hold, so $T$ is positively expansive.
	\end{proof}
	\begin{remark}\label{diaexp} \rm 
		%(a) The diagonal operator $A$ in Corollary \ref{diagcor} is bounded if and only if
		%		\[
		%		\sup_{n\ge1}|\lambda_n|<\infty.
		%		\]
		%		Indeed, if $C:=\sup_n|\lambda_n|<\infty$, then
		%		\[
		%		\|Ax\|\le C\|x\|
		%		\]
		%		on either $\ell^p$ or $c_0$. Conversely, since the standard unit vectors belong to $D(A)$,
		%		\[
		%		\|Ae_n\|=|\lambda_n|,\quad \|e_n\|=1.
		%		\]
		%		Thus boundedness of $A$ requires the sequence $(|\lambda_n|)$ to be bounded.
		
		(a) The diagonal semigroup in Corollary \ref{diagcor} is uniformly exponentially expanding if and only if
		\[
		\inf_{n\ge1}\re \lambda_n>0.\]
		Indeed, put $\alpha=\inf_n\re \lambda_n$. For the diagonal semigroup,
		\[
		\inf_{\|x\|=1}\|T(t)x\|
		=\inf_n e^{t\re \lambda_n}.
		\]
		The lower estimate follows coordinatewise, and the reverse inequality follows by taking $x=e_n$.
		If $\alpha>0$, then
		\[
		\|T(t)x\|\ge e^{\alpha t}\|x\|,
		\]
		so $T$ is uniformly exponentially expanding.
		Conversely, suppose there exist $M,a>0$ such that
		\[
		\|T(t)x\|\ge Me^{at}\|x\|
		\quad \mbox{for}\quad t\ge0.
		\]
		Taking $x=e_n$ gives
		\[
		e^{t\re \lambda_n}\ge Me^{at}.
		\]
		For $t>0$,
		\[
		\re \lambda_n\ge a+\frac{\log M}{t}.
		\]
		Letting $t\to\infty$, we obtain $\re \lambda_n\ge a$ for every $n$. Hence $\alpha\ge a>0$.
		
		(b)  Choose $\lambda_n=\frac1n+i$ in Corollary \ref{diagcor}, we have 
		$\re\lambda_n=\frac1n>0$ for all $n$ and $\displaystyle\inf_{n}|\lambda_n|=\displaystyle\inf_n\sqrt{1+\frac{1}{n^2}}=1$. So, the associated diagonal semigroup is positively expansive, but it is not uniformly exponentially expanding. This shows that a positively expansive semigroup does not need to be uniformly exponentially expanding.
		
		(c) Choose $\lambda_n=\frac1n$ in Corollary \ref{diagcor}. Then $\re\lambda_n>0$ for all $n$, but $\inf_n|\lambda_n|=0$, so the associated diagonal semigroup is positively kinematic expansive but not positively expansive. Together with (b), this shows that both implications in \eqref{3notion}
		are strict.		
	\end{remark}

	\section{Spectral criteria for the generator}\label{sec6} 
	In this section we give criteria for expansiveness in terms of the generator. For general generators, we compare the lower growth bound with the approximate point spectrum of the generator and show that they coincide under an approximate spectral mapping hypothesis. For bounded generators, we characterize uniform exponential expansion by the location of the spectrum and give a necessary spectral condition for positive expansiveness.
	\subsection{General generators} 
	To compare the lower growth bound of $T$ with the spectrum of $A$, we use the following quantity. For the generator $A$ of a $C_0$-semigroup, put
	\begin{equation*}
		{\rm i}_\ap(A):=\inf\{\re\lambda:\ \lambda\in\sigma_\ap(A)\} 
	\end{equation*}
	with convention $\inf\varnothing:=+\infty $.

	\begin{theorem}\label{geAthm}
		Let $T$ be a $C_0$-semigroup with generator $A$. The following holds:
		\begin{enumerate}
			\item[(i)] For every $t>0$, $e^{t\sigma_\ap(A)}\subset\sigma_\ap(T(t))$. Consequently, $\omega_-(T)\le {\rm i}_\ap(A)$. In particular, if \eqref{mthmeq} holds for constants $M,a>0$, then ${\rm i}_\ap(A)\ge a$.
			\item[(ii)] Suppose that there exists $\tau>0$ such that
			\begin{equation}\label{geAeq1}
				\sigma_\ap(T(\tau))\subset e^{\tau\sigma_\ap(A)} .
			\end{equation}
			Then $\omega_-(T)={\rm i}_\ap(A)$. Consequently, $T$ is uniformly exponentially expanding if and only if ${\rm i}_\ap(A)>0$, and in that case $T$ is positively expansive.
		\end{enumerate}
	\end{theorem}
	
	\begin{proof}
		(i) The first statement is the spectral inclusion theorem; see \cite[IV.3.6]{EN}. Next, if $\sigma_\ap(A)=\varnothing$, the inequality  $\omega_-(T)\leq {\rm i}_\ap(A)$ is trivial. If $\lambda\in\sigma_\ap(A)$, then  by \eqref{Lmeq}, \[\tn{T(t)}\le|e^{t\lambda}|=e^{t\re\lambda}\quad\mbox{ for}\quad t>0,\] that is, 
		\[\frac{1}{t}\log\tn{T(t)}\le\re\lambda\quad\mbox{ for}\quad t>0.\] 
		Taking the supremum over $t>0$ and using \eqref{omega-T} gives $\omega_-(T)\le\re\lambda$. As $\lambda\in\sigma_\ap(A)$ was arbitrary, $$\omega_-(T)\le \inf\{\re\lambda:\ \lambda\in\sigma_\ap(A)\}={\rm i}_\ap(A).$$ 
		
		 Finally, if \eqref{mthmeq} holds for some $M,a>0$, then \[\tn{T(t)}\ge Me^{at}\quad\mbox{for}\quad t\geq 0\] and $\omega_-(T)\ge a$ by the definition of $\omega_-(T)$, so ${\rm i}_\ap(A)\geq a$. 
		
		(ii) Since $T(\tau)$ is bounded, we have $\sigma_\ap(T(\tau))\ne\varnothing$, so \eqref{geAeq1} implies $\sigma_\ap(A)\neq\varnothing$. Every $\mu\in\sigma_\ap(T(\tau))$ has the form $\mu=e^{\tau\lambda}$ with $\lambda\in\sigma_\ap(A)$, so $|\mu|=e^{\tau\re\lambda}\ge e^{\tau {\rm i}_\ap(A)}$ (with $e^{-\infty}:=0$). By Lemma \ref{lgblm}\,(iii) and \eqref{MZeq},
		\begin{equation*}
			e^{\tau\omega_-(T)}=r_-(T(\tau))=\min_{\mu\in\sigma_\ap(T(\tau))}|\mu|\ \ge\ e^{\tau {\rm i}_\ap(A)} ,
		\end{equation*}
		so $\omega_-(T)\ge {\rm i}_\ap(A)$. Together with (i), $\omega_-(T)={\rm i}_\ap(A)$. The last assertions follow from Theorem \ref{uepeqthm}\,(iii).
	\end{proof}
	
	\begin{remark}\rm 
		By Theorem \ref{geAthm}\,(i), the inclusion $e^{\tau\sigma_\ap(A)}\subset\sigma_\ap(T(\tau))$ always holds, so \eqref{geAeq1} is equivalent to the equality $\sigma_\ap(T(\tau))=e^{\tau\sigma_\ap(A)}$. Since $\sigma_\ap(T(\tau))$ is closed, so is $e^{\tau\sigma_\ap(A)}$, and  $T(\tau)$ is bounded from below since $0\not\in e^{\tau\sigma_\ap(A)}$, i.e. $\tn{T(\tau)}>0$. Lemma \ref{lgblm}\,(i) implies that $\tn{T(t)}>0$ for all $t\geq 0$. This holds if $A$ is bounded, and for normal generators on Hilbert spaces and for diagonal semigroups precisely when $e^{\tau\sigma_\ap(A)}$ is closed. It may fail even for diagonal semigroups: for the entries $\lambda_n=\frac1n+\frac{2\pi i n}{\tau}$ one has $1\in\sigma_\ap(T(\tau))\setminus e^{\tau\sigma_\ap(A)}$; see Corollary \ref{diagcor} for the definition of diagonal semigroup.
		{\hfill$\Diamond$}
	\end{remark}
	 
\subsection{Bounded generators}

	For uniformly continuous semigroups the conditions of Proposition \ref{T-1pro} can be expressed in form of $\sigma(A)$. Recall that a semigroup $(T(t))_{t\ge 0}$ is uniformly continuous if and only if its generator $A$ is bounded and $T(t)=e^{tA}, t\ge 0$; see \cite[Theorem~II.3.7]{EN}. 
	
	\begin{theorem}\label{Abdeqthm}
		Let $A\in\LL(X)$ and $T(t)=e^{tA},\, t\ge 0$. The following assertions are equivalent:
		\begin{enumerate}
			\item[(i)] $\sigma(A)\subset\{\lambda\in\C:\re\lambda>0\}$;
			\item[(ii)] $\alpha:=\min\{\re\lambda:\lambda\in\sigma(A)\}>0$, equivalently $\displaystyle\sup\{\re\mu:\mu\in\sigma(-A)\}<0$;
			\item[(iii)] $T$ satisfies the equivalent conditions of Proposition \ref{T-1pro}.
		\end{enumerate}
		In this case $T$ extends to a uniformly continuous group $(e^{tA})_{t\in\R}$, $T$ is uniformly exponentially expanding, hence positively expansive, and $r(e^{-A})=e^{-\alpha}$.
	\end{theorem}
	
	\begin{proof} (i)~$\Leftrightarrow$~(ii). 	Since $A$ is bounded, its spectrum  $\sigma(A)$ is nonempty and compact, so $\re$ attains its minimum on it; this gives (i)~$\Leftrightarrow$~(ii), and the equivalent form follows from $\sigma(-A)=-\sigma(A)$.
		
		(ii)~$\Rightarrow$~(iii). The exponential $e^{tA}$ is defined for all $t\in\R$ and $e^{tA}$ is invertible with inverse $e^{-tA}$. By the spectral mapping theorem for uniformly continuous semigroup of bounded operators (see \cite[Theorem I.3.13]{EN}),  \[\sigma(e^{-A})=\{e^{-\lambda}:\lambda\in\sigma(A)\},\]
		 so \[r(e^{-A})=\max_{\lambda\in\sigma(A)}e^{-\re\lambda}=e^{-\alpha}<1.\] By Gelfand's formula there is $n_0\ge1$ with \[\|T(n_0)^{-1}\|=\|e^{-n_0A}\|=\|(e^{-A})^{n_0}\|<1,\] which is Proposition \ref{T-1pro}\,(i).
		
		(iii)~$\Rightarrow$~(i). By Proposition \ref{T-1pro}\,(v) there is $\tau>0$ with $\sigma(e^{\tau A})\subset\{|\mu|>1\}$. Since $\sigma(e^{\tau A})=e^{\tau\sigma(A)}$, every $\lambda\in\sigma(A)$ satisfies $e^{\tau\re\lambda}>1$, i.e.\ $\re\lambda>0$.
		
		The final statements follow from Proposition \ref{T-1pro} and the computation of $r(e^{-A})$.
	\end{proof}
	
	\begin{example}[{\cite[Example 4]{LM}}]\label{scaex}\rm 
		On any nonzero complex Banach space $X$, consider semigroup $T(t)x=2^tx$ for $x\in X, t\geq 0$. Its generator is $A=(\log 2)I\in \LL(X)$, so $T$ is uniformly continuous. Its spectrum $\sigma(A)=\{\log2\}$ and we can apply Theorem \ref{Abdeqthm} to obtain positive expansiveness of $T$. Besides,  $\|T(1)^{-1}\|=\tfrac12<1$, so Proposition \ref{T-1pro} also applies. We do not need to calculate the norm of resolvent  as in \cite{LM}.
	\end{example}
	
	\begin{corollary}[Finite dimensions]\label{fdcor}
		Let $X$ be a finite-dimensional complex Banach space and $A\in\LL(X)$. The following assertions are equivalent:
		\begin{enumerate}
			\item[(i)] $e^{tA}$ is positively kinematic expansive;
			\item[(ii)] $e^{tA}$ is positively expansive;
			\item[(iii)] $\sigma(A)\subset\{\lambda\in\C:\re\lambda>0\}$;
			\item[(iv)] $e^{tA}$ is uniformly exponentially expanding.
		\end{enumerate}
	\end{corollary}
	
	\begin{proof}
		(i)~$\Rightarrow$~(iii). In finite dimensions $\sigma(A)=\sigma_p(A)$. Assume that $e^{tA}$ is positively kinematic expansive. Theorem \ref{nethm}\,(iii) gives $\sigma_p(A)\subset\{\re\lambda>0\}$ and hence (iii) holds. 
		
		(iii)~$\Rightarrow$~(iv).
		 This follows from Theorem \ref{Abdeqthm} together with Proposition \ref{T-1pro}. 
		 
		 (iv)~$\Rightarrow$~(ii). This is Theorem \ref{mainthm}.
		 
		 (ii)~$\Rightarrow$~(i). This is obvious. 		 
	\end{proof}

	In infinite dimensions the spectrum of a bounded generator might not consist of eigenvalues, and positive expansiveness yields a weaker spectral restriction than uniform exponential expansion: the open right half-plane of Theorem~\ref{Abdeqthm} is replaced by the closed one with the origin removed.
	\begin{theorem}\label{Abdthm}
		Let $A\in\LL(X)$ on a complex Banach space and suppose $T(t)=e^{tA}$ is positively expansive. Then $A$ is invertible and $$\sigma(A)\subset\{\lambda\in\C:\re\lambda\ge 0\}\setminus\{0\}.$$ Consequently, $$\omega_-(T)=\min\{\re\lambda:\lambda\in\sigma(A)\}\ge 0.$$
	\end{theorem}
	In order to prove this theorem, we need the following result.
	\begin{lemma}\label{zerorealthm}
		Let $T$ be a $C_0$-semigroup with generator $A$. If $T$ is positively expansive, then no point $i\beta$ with $\beta\in\R$ is a limit of points of $\sigma_\ap(A)\cap\{\re\lambda<0\}$, i.e.
		\begin{equation}
			\overline{\sigma_{\mathrm{ap}}(A)\cap
					\{\lambda\in\mathbb C:\re\lambda<0\}}
				\subseteq
				\{\lambda\in\mathbb C:\re\lambda<0\}.
		\end{equation}
	\end{lemma}
	\begin{proof} 
	Suppose $\lambda_k=-\alpha_k+i\beta_k\in\sigma_\ap(A)$ with $\alpha_k>0$ and $\lambda_k\to i\beta$. Since $\sigma_\ap(A)$ is closed, $i\beta\in\sigma_\ap(A)$, so $\beta\ne0$ by Theorem \ref{nAxthm}. Discarding finitely many $k$ we may assume \[\frac{|\beta|}2\le|\beta_k|\le|\beta|+1\quad\mbox{and}\quad \alpha_k\to 0.\]
	Put $h_k=\frac{2\pi}{|\beta_k|}$ and $\eps=\frac{2\pi}{|\beta|+1}$, so that $h_k\ge\eps$ and $e^{i\beta_kh_k}=1$. Let $\delta\in(0,1)$ correspond to $\eps$ in Definition \ref{LMdef}. Fix $k$ so large that $\alpha_kh_k<\delta/4$, and write $\lambda=\lambda_k$, $\alpha=\alpha_k$, $h=h_k$. Choose $t_0>0$ with $2e^{-\alpha t_0}<\delta/8$, and let
	\begin{equation*}
		C:=\sup_{0\le u\le t_0+1+h}\int_0^ue^{(u-r)\re\lambda}\|T(r)\|\,dr<\infty .
	\end{equation*}
	Since $\lambda\in\sigma_\ap(A)$, we can choose unit vectors $z_n\in D(A)$ with $(A-\lambda I)z_n\to 0$ as $n\to\infty$. Fix $u>0$. Since $z_n\in D(A)$, the orbit $r\mapsto T(r)z_n$ is continuously differentiable with derivative $T(r)Az_n$. By the product rule $\phi_n(r):=e^{\lambda(u-r)}T(r)z_n$ is continuously differentiable on $[0,u]$ with
	\begin{equation*}
		\phi_n'(r)=-\lambda e^{\lambda(u-r)}T(r)z_n+e^{\lambda(u-r)}T(r)Az_n=e^{\lambda(u-r)}T(r)(A-\lambda I)z_n .
	\end{equation*}
	Integrating $\phi_n(u)-\phi_n(0)=\int_0^u\phi_n'(r)\,dr$ gives
	\begin{equation*}
		T(u)z_n-e^{u\lambda}z_n=\int_0^ue^{\lambda(u-r)}T(r)(A-\lambda I)z_n\,dr,
	\end{equation*}
	so, for $0\le u\le t_0+1+h$
	\begin{align}\label{zerorealeq}
		\|T(u)z_n-e^{u\lambda}z_n\|
		\le \left(\int_0^ue^{(u-r)\re\lambda}\|T(r)\|dr\right)\|(A-\lambda I)z_n\|
		\leq C\theta_n,
	\end{align} 
	where \[ \theta_n:=\|(A-\lambda I)z_n\|\to 0 \quad\mbox{as}\quad  n\to\infty.\]
	
	Fix $n$ and put $x=T(h)z_n$, $y=z_n$, and define $s\in{\cal H}$ by
	\[s(t)=\begin{cases}
		t&\mbox{for}\quad 0\le t\le t_0,\\
		t+h(t-t_0)&\mbox{for}\quad  t_0\le t\le t_0+1,\\
		t+h&\mbox{for}\quad t\ge t_0+1.			 
	\end{cases} \]  We estimate \[T(t)x-T(s(t))y=T(t+h)z_n-T(s(t))z_n.\]
	For $t\ge t_0+1$ we have $s(t)=t+h$, so $T(t)x-T(s(t))y=0$.
	
	For $0\le t\le t_0$ we have $s(t)=t$, and $t,\,t+h\in[0,t_0+1+h]$. By \eqref{zerorealeq} and the triangle inequality,
	\begin{align*}
		\|T(t+h)z_n-T(t)z_n\|
		&\le\|T(t+h)z_n-e^{\lambda(t+h)}z_n\|+\|e^{\lambda(t+h)}z_n-e^{\lambda t}z_n\|\\
		&+\|e^{\lambda t}z_n-T(t)z_n\|
		\le |e^{\lambda t}|\,|e^{\lambda h}-1|+2C\theta_n .
	\end{align*}
	Since $e^{i\beta_kh}=1$, we have $e^{\lambda h}=e^{-\alpha h}$, hence
	\[
	|e^{\lambda t}|\,|e^{\lambda h}-1|=e^{-\alpha t}\,(1-e^{-\alpha h})\le\alpha h<\frac\delta4 ,
	\]
	and therefore \[\|T(t)x-T(s(t))y\|\le\frac\delta4+2C\theta_n.\]
	
	For $t_0\le t\le t_0+1$ we write $s(t)=t+d$ with $d=h(t-t_0)\in[0,h]$; then $t+d,\,t+h\in[0,t_0+1+h]$. Again by \eqref{zerorealeq},
	\begin{align*}
		\|T(t+h)z_n-T(t+d)z_n\|
		&\le|e^{\lambda(t+h)}|+|e^{\lambda(t+d)}|+2C\theta_n\\
		&=e^{-\alpha(t+h)}+e^{-\alpha(t+d)}+2C\theta_n\\
		&\le\ 2e^{-\alpha t_0}+2C\theta_n<\frac\delta8+2C\theta_n .
	\end{align*}
	
	Consequently,
	\[
	\|T(t)x-T(s(t))y\|\le\frac\delta4+2C\theta_n\quad\mbox{for all}\quad t\ge0 .
	\]
	Since $C$ does not depend on $n$ and $\theta_n\to0$ as $n\to\infty$, we may fix $n$ with $C\theta_n<\delta/8$. Then \[\|T(t)x-T(s(t))y\|\le\frac\delta4+2C\theta_n<\frac\delta2\quad\mbox{for all}\quad t\ge0.\] 
	Moreover, since $\alpha h<\frac\delta4<\frac14$ and $C\theta_n<\frac18$,
	\[
	\|x\|=\|T(h)z_n\|\ge e^{-\alpha h}-C\theta_n\ge1-\alpha h-C\theta_n>\tfrac12 .
	\]
	Hence
	\[ \|T(t)x-T(s(t))y\|<\delta\|x\|\quad\mbox{for all}\quad t\ge 0,\]
	and positive expansiveness yields $t_*\in[0,\eps)$ with $y=T(t_*)x$ or $x=T(t_*)y$. In the first case $z_n=T(t_*+h)z_n$ with $t_*+h>0$, so $z_n$ is a nonzero periodic point, contradicting Theorem \ref{nethm}\,(ii). In the second case $T(h)z_n=T(t_*)z_n$, and injectivity (Theorem \ref{nethm}\,(i)) gives $T(h-t_*)z_n=z_n$ with $h-t_*>0$ (since $t_*<\eps\le h$), again a nonzero periodic point. This contradiction proves the lemma.
	\end{proof}

	\begin{proof}[Proof of Theorem \ref{Abdthm}] Suppose on the contrary that $$\sigma(A)\cap\{\lambda\in\C:\re\lambda<0\}\ne\varnothing.$$
We first show that  $$\sigma_\ap(A)\cap \{\lambda\in\C:\re\lambda<0\}\ne\varnothing.$$ 
 Since $\sigma(A)$ is compact, it contains a point $\lambda_0$ of minimal real part, and $\re\lambda_0<0$. Every neighbourhood of $\lambda_0$ contains points $\lambda_0-\eta$, $\eta>0$, with  $\lambda_0-\eta\not\in\sigma(A)$; hence $\lambda_0\in\partial\sigma(A)\subset\sigma_\ap(A)$, recalling \eqref{sigmaL}. Thus
\[
\Sigma:=\sigma_\ap(A)\cap\{\lambda\in\C:\re\lambda<0\}
\]
is nonempty, and it is bounded because $\sigma_\ap(A)\subset\sigma(A)$.

 By Lemma~\ref{zerorealthm}, $\overline{\Sigma}$ is a compact subset of $\{ \lambda\in\C: \re\lambda<0\}$. Consequently, there exists
			$\gamma>0$ such that
			\[
			\Sigma\subset\{\lambda\in\C:\re\lambda\le-2\gamma\},\quad\mbox{i.e.,}\quad \sigma_\ap(A)\cap\{\lambda\in\C:-2\gamma<\re\lambda<0\}=\varnothing.
			\]
			The line $\ell=\{-\gamma+is:s\in\R\}$ therefore avoids
			$\sigma_\ap(A)$ and, in particular, $\partial\sigma(A)$.
			It follows that $\ell$ is the disjoint union of the relatively
			open sets
			\[
			\ell\cap\operatorname{int}\sigma(A)
			\quad\text{and}\quad
			\ell\cap\rho(A),
			\]
			both relatively open in $\ell$. The second set is nonempty, since $\lambda\in\rho(A)$ whenever $|\lambda|>\|A\|$. As $L$ is connected, $\ell\cap\operatorname{int}\sigma(A)=\varnothing$, so $\ell\subset\rho(A)$.
			Thus $\sigma(A)$ splits into the disjoint compact sets\[
			\sigma_1:=\sigma(A)\cap\{\lambda\in\C:\re\lambda<-\gamma\},\quad \sigma_2:=\sigma(A)\cap\{\lambda\in\C:\re\lambda>-\gamma\},
			\]
			and $\sigma_1\supset\Sigma\ne\varnothing$. Let $P\ne0$ be the Riesz projection of $A$ associated with $\sigma_1$. The restriction
			$A_1:=A|_{\ran P}$ satisfies $\sigma(A_1)=\sigma_1$; see
			\cite[Section~III.6.4]{Kato}. Since
			$$\max\{\lambda\in\C:\re\sigma(A_1)\}<-\gamma,$$ the equality of the spectral
			and growth bounds for bounded generators
			\cite[Corollary~IV.2.4]{EN} gives
			\[
			\|e^{tA_1}\|\le Ce^{-\gamma t}
			\quad\mbox{for}\quad t\ge0,
			\]
			for some $C\ge1$. Every nonzero $x\in\ran P$ then has a
			bounded orbit under $T(t)=e^{tA}$, contradicting
			Theorem~\ref{ubothm}.
			 Therefore
			\[
			\sigma(A)\subset\{\lambda\in\C:\re\lambda\ge0\}.
			\]
			
			It remains to show that $0\not\in\sigma(A)$. 
			If $0\in\sigma(A)$, the preceding inclusion implies
			$0\in\partial\sigma(A)\subset\sigma_\ap(A)$, contrary to
			Theorem~\ref{nAxthm}, since $T$ is positively expansive. Hence $A$ is invertible and
			\[
			\sigma(A)\subset\{\lambda\in\C:\re\lambda\ge0\}\setminus\{0\}.
			\]
			
			Finally, $T(t)^{-1}=e^{-tA}$ and \eqref{L-1} give
			\[
			\tn{T(t)}=\|e^{-tA}\|^{-1}.
			\]
			Since $-A$ is bounded, its semigroup has growth bound
			$s(-A)$ by \cite[Corollary~IV.2.4]{EN}. Consequently,
			\[
			\begin{aligned}
				\omega_-(T)
				&=-\lim_{t\to\infty}\frac1t\log\|e^{-tA}\|\\
				&=-s(-A)
				=\min_{\lambda\in\sigma(A)}\re\lambda
				\ge0.
			\end{aligned}
			\]
			This completes the proof.
		
	\end{proof}

The gap between Theorems~\ref{Abdeqthm} and~\ref{Abdthm} is genuine, and the necessary conditions of Section~\ref{sec3} are not sufficient for positive expansiveness. The following semigroup is positively kinematic expansive and has empty point spectrum, but it is not positively expansive, because $0\in\sigma(A)$ contradicts Theorem~\ref{Abdthm}.
\begin{example}\label{mulex}
\rm Let $X=C([0,1];\mathbb C)$, equipped with the supremum norm, and define an operator $A$ on $X$ by
	\[
	(Af)(x)=xf(x),\ f\in X,\ x\in [0,1].
	\]
	Then  $A$ is bounded, with $\|A\|=1$, and generates the uniformly continuous semigroup $T$ given by
	\[
	(T(t)f)(x)=e^{tx}f(x),\ f\in X, \, x\in [0,1],\, t\ge 0.
	\]
	
	(a) $T$ satisfies all three necessary conditions in Theorem \ref{nethm}. Indeed, if $T(t)(f)=0$ for some $t>0$, then $f(x)=0$ for all $x\in [0,1]$. It follows that $T(t)$ is injective for all $t>0$.  If $T(t)f=f$ for some $t>0$, then
	$(e^{tx}-1)f(x)=0$.
	Hence $f(x)=0$ for $x\in (0,1]$, and the continuity gives $f(0)=0$, so $f(x)=0$ for all $x\in [0,1]$. Thus, $T$ has no nonzero periodic points. 
	In addition, the point spectrum
	$\sigma_p(A)=\varnothing$ since if $(\lambda I-A)f=0$, then $(x-\lambda)f(x)=0$ for all $x\in [0,1]$. The continuity of $f$ implies that $f(x)=0$ for all $x\in [0,1]$.	
	
	(b) Every nonzero orbit is unbounded. For, if $f\ne0$, then $f(x_0)\ne0$ for some $x_0>0$. Consequently,
	\[
	\|T(t)f\|_\infty
	\ge e^{tx_0}|f(x_0)|
	\to\infty.
	\]
	
	(c) $T$ is not positively expansive by Theorem \ref{Abdthm} since $\sigma(A)=[0,1]$. Indeed, if $\lambda\in [0,1]$, then
	$	((\lambda I-A)f)(\lambda)=0$ for all $f\in X$.
	If $\lambda\notin[0,1]$, then
	$	((\lambda I-A)^{-1}g)(x)=\frac{g(x)}{\lambda-x}$ for  $g\in X, 
	x\in [0,1]. 
	$
	Thus,  $\sigma(A)=[0,1]$.
\end{example}

	\section*{Acknowledgments}
	This work was carried out during the author's stay at the Vietnam Institute for Advanced Study in Mathematics (VIASM). The author thanks VIASM for its hospitality and excellent working conditions.

%\tableofcontents

\begin{thebibliography}{99}
	
	\bibitem{Art} A. Artigue, \textit{Kinematic expansive flows}, Ergodic Theory Dynam.  Systems \textbf{36} (2016), 390--421.
	
	\bibitem{BCDMP} N.C. Bernardes Jr., P.R. Cirilo, U.B. Darji, A. Messaoudi and E.R. Pujals, \textit{Expansivity and shadowing in linear dynamics}, J. Math. Anal. Appl. \textbf{461} (2018), 796--816.
	
	\bibitem{BM} N.C. Bernardes Jr.  and A. Messaoudi, \textit{Shadowing and structural stability for operators}, Ergodic Theory Dynam.  Systems \textbf{41} (2021), 961--980.
	
	\bibitem{BW} R. Bowen and P. Walters, \textit{Expansive one-parameter flows}, J.  Differential Equations \textbf{12} (1972), 180--193.
	\bibitem{CGP} P.R. Cirilo, B. Gollobit and E. R. Pujals, \textit{Dynamics of generalized hyperbolic linear operators}, Adv.  Math.  \textbf{387} (2021), Paper No.~107830.
	
	\bibitem{Conway} 	J.B. Conway,
	\textit{A Course in Functional Analysis},
	2nd ed., Graduate Texts in Mathematics, vol.~96,
	Springer-Verlag, New York, 1990.
	
%	\bibitem{Datko} R. Datko, \textit{Extending a theorem of A. M. Liapunov to Hilbert space}, J.  Math.  Anal.  Appl.  \textbf{32} (1970), 610--616.
	
	%\bibitem{DSW} W. Desch, W. Schappacher and G. F. Webb, \textit{Hypercyclic and chaotic semigroups of linear operators}, Ergodic Theory Dynam.  Systems \textbf{17} (1997), 793--819.
	
	\bibitem{Eis} M. Eisenberg, \textit{Expansive automorphisms of finite dimensional vector spaces}, Fund.  Math.  \textbf{59} (1966), 307--312.
	
	\bibitem{EH} M. Eisenberg and J. H. Hedlund, \textit{Expansive automorphisms of Banach spaces}, Pacific J.  Math.  \textbf{34} (1970), 647--656.
	
	\bibitem{EN} K.-J. Engel and R. Nagel, \textit{One-Parameter Semigroups for Linear Evolution Equations}, Graduate Texts in Mathematics, vol. 194, Springer-Verlag, New York, 2000.
	
	\bibitem{GEM} K.-G. Grosse-Erdmann and A.P. Manguillot, \textit{Linear Chaos}, Universitext, Springer, London, 2011.
	
	\bibitem{Hed} J.H. Hedlund, \textit{Expansive automorphisms of Banach spaces, II}, Pacific J.  Math.  \textbf{36} (1971), 671--675.
	
	\bibitem{HJ} S. Herrero and N. Jaque, \textit{On two notions of expansiveness for continuous semiflows}, J. Math. Anal. Appl. \textbf{515(1)}(2022), 126405.
	
	\bibitem{Hien} H.M. Huynh, \emph{Katok--Hasselblatt-kinematic expansive flows}, J. Korean Math. Soc. \textbf{59} (2022), no.~1, 151--170.
	
	%\bibitem{HP} E. Hille and R.S. Phillips, \textit{Functional Analysis and Semi-Groups}, Amer. Math. Soc. Colloq. Publ. \textbf{31}, Amer. Math. Soc., Providence, RI, 1957.	
	\bibitem{JNY} W. Jung, N. Nguyen and Y. Yang, \textit{Spectral decomposition for rescaling expansive flows with rescaled shadowing}, Discrete Contin.  Dyn.  Syst.  \textbf{40} (2020), 2267--2283.
	\bibitem{Kato} T. Kato, \emph{Perturbation Theory for Linear Operators}, reprint of the 1980 edition, Classics in Mathematics, Springer, Berlin, 1995.
	\bibitem{LM} J. Lee and C.A. Morales, \textit{Bowen--Walters expansiveness for semigroups of linear operators}, Ergodic Theory Dynam.  Systems \textbf{43} (2023), 1942--1951.
	
	\bibitem{LMN} K. Lee, C.A. Morales and N. Nguyen, \textit{Topologically stable linear operators}, Bull.  Soc.  Math.  France \textbf{151} (2023), 647--657.
	
	\bibitem{LMN2} K. Lee, C.A. Morales and N. Nguyen, \textit{Topological Dynamics from a Measure-Theoretic Viewpoint}, Springer Asia Pacific Mathematics Series, vol.~14, Springer, Singapore, 2026.
	
	\bibitem{LMV} K. Lee, C.A. Morales and H. Villavicencio, \textit{Asymptotic expansivity}, J.  Math.  Anal.  Appl.  \textbf{507} (2022), Paper No.~125729.
	
	\bibitem{LR} K. Lee and A. Rojas, \textit{Eventually expansive semiflows}, Commun.  Pure Appl.  Anal.  \textbf{21} (2022), 3335--3351.
	
	\bibitem{Mai} M. Maiuriello, \textit{Expansivity and strong structural stability for composition operators on $L^p$ spaces}, Banach J.  Math.  Anal.  \textbf{16} (2022), Paper No.~51.
	
	\bibitem{MZ} E. Makai and J. Zem\'anek, \textit{The surjectivity radius, packing numbers and boundedness below of linear operators}, Integral Equations Operator Theory \textbf{6} (1983), 372--384.
	
	%\bibitem{Mul} V. M\"uller, \textit{Spectral Theory of Linear Operators and Spectral Systems in Banach Algebras}, 2nd ed., Operator Theory: Advances and Applications, vol.~139, Birkh\"auser, Basel, 2007.
	
	%\bibitem{vN} J. van Neerven, \textit{The Asymptotic Behaviour of Semigroups of Linear Operators}, Operator Theory: Advances and Applications, vol.~88, Birkh\"auser, Basel, 1996.
	
	\bibitem{Pazy} A. Pazy, \textit{Semigroups of Linear Operators and Applications to Partial Differential Equations}, Applied Mathematical Sciences, vol. 44, Springer-Verlag, New York, 1983.
	
	\bibitem{Utz} W.R. Utz, \textit{Unstable homeomorphisms}, Proc.  Amer.  Math.  Soc.  \textbf{1} (1950), 769--774.
	
	\bibitem{WW} X. Wen and L. Wen, \textit{A rescaled expansiveness of flows}, Trans.  Amer.  Math.  Soc.  \textbf{371} (2019), 3179--3207.
	
	\bibitem{YSMN} Y. Yang, J. Shen, C.A. Morales, N. Nguyen,
	\textit{Expansivity for Cauchy problems on Banach spaces},
	J. Differential Equations,
	\textbf{471}(2026), Paper No. 114403
	
\end{thebibliography}
\end{document}